\documentclass[11pt,reqno]{amsart}
\usepackage{amssymb,amsmath,amsthm}

\newtheorem{thm}{Theorem}%[section]
\newtheorem{prop}[thm]{Proposition}
\newtheorem{lem}[thm]{Lemma}
\newtheorem{cor}[thm]{Corollary}

\newtheorem{defn}[thm]{Definition}
\newtheorem{rem}[thm]{Remark}
\newcommand{\n}{\noindent}
\newcommand{\va}{\vskip1pc}
\newcommand{\R}{\mathbb R}

\newcommand{\ZZ}{{\bar Z}}

\newcommand{\C}{\mathbb C}

\DeclareMathOperator{\Imm}{Im}
\DeclareMathOperator{\supp}{supp}

\DeclareMathOperator{\Rre}{Re}
\DeclareMathOperator{\Dom}{Dom}
\DeclareMathOperator{\Ran}{Range}

\newcommand{\p}{\partial}

\newcommand{\les}{\lesssim}
\newcommand{\ges}{\gtrsim}
\newcommand{\z}{\bar z}
\newcommand{\w}{\bar w}
\newcommand{\dbar}{\bar\partial}
\newcommand{\vp}{\varphi}

\newcommand{\atopp}[2]{\genfrac{}{}{0pt}{2}{#1}{#2}}

\newcommand{\A}[2]{a_{#1}^{#2}}
\newcommand{\sjk}{\sum_{j,k\geq 1}}

\newcommand{\nn}{\nonumber}

\newcommand{\Ts}{T_\tau}
\newcommand{\Ks}{K_\tau}

\newcommand{\Zstp}{Z_{\tau p}}
\newcommand{\Zbstp}{\bar{Z}_{\tau p}}

\newcommand{\Zstpz}{Z_{\tau p,z}}
\newcommand{\Zbstpz}{\bar{Z}_{\tau p,z}}
\newcommand{\Ms}{M_{\tau p}}

\newcommand{\Wstp}{W_{\tau p}}
\newcommand{\Wbstp}{\overline{W}_{\tau p}}
\newcommand{\Wstpw}{W_{\tau p,w}}
\newcommand{\Wbstpw}{\overline{W}_{\tau p,w}}

\newcommand{\Kse}{K_{\tau,\epsilon}}

\newcommand{\Ystp}{Y_{\tau p}}
\newcommand{\Xstp}{X_{\tau p}}
\newcommand{\epl}{e^{i \tau T(w,z)}}
\newcommand{\emi}{e^{-i \tau T(w,z)}}
\newcommand{\T}{\epl \frac{\p}{\p\tau} \emi}

\newcommand{\ep}{\epsilon}
\newcommand{\sh}{^{\#}}
\newcommand{\G}{\tilde G}

\newcommand{\Boxtp}{\Box_{\tau p}}

\newcommand{\Boxtpz}{\Box_{\tau p,z}}
\newcommand{\Boxtpw}{\Box_{\tau p,w}}
\newcommand{\Boxtpx}{\Box_{\tau p,x}}
\newcommand{\Boxtpy}{\Box_{\tau p,y}}

\newcommand{\diam}[1]{\text{diam}(#1)}

\newcommand{\HH}{\mathcal{H}}

\renewcommand{\H}{\mathcal{H}_{\tau p}}
\newcommand{\Hw}{\tilde{\mathcal{H}}_{\tau p}}
\newcommand{\Htp}{H_{\tau p}}
\newcommand{\Hwtp}{{\tilde H}_{\tau p}}

\newcommand{\lam}{\lambda}

\newcommand{\toap}{\textit{To appear}, J. Funct. Anal.}
 
\begin{document}

\title{Heat Equations in  $\R\times\C$\footnote{\toap}}

\author{Andrew S. Raich}

\address{Department of Mathematics, Texas A\&M University, Mailstop 3368, 
College Station, TX 77843-3368}

\begin{abstract}
Let $p:\C\to\R$ be a subharmonic, nonharmonic polynomial and $\tau\in\R$ a parameter.
Define $\Zbstp = \frac{\p}{\p\z} + \tau\frac{\p p}{\p\z}$, a closed, densely-defined operator
on $L^2(\C)$. If  $\Boxtp = \Zbstp\Zbstp^*$ and $\tau>0$,
we solve the heat equation
$\frac{\p u}{\p s} +  \Boxtp u =0$, $u(0,z) = f(z)$,
on $(0,\infty)\times\C$. The solution 
comes via the heat semigroup $e^{-s\Boxtp}$, and we show
that $u(s,z) = e^{-s\Boxtp}[f](z) = \int_{\C}\Htp(s,z,w) f(w)\, dw$.
We prove that $\Htp$ is $C^\infty$ off the diagonal $\{(s,z,w):s=0 \text{ and }z=w\}$
and that $\Htp$ and its derivatives have exponential decay. In particular,
we give new estimates for the long time behavior of the heat equation.
\end{abstract}

\keywords{
heat kernel, weighted $\bar\partial$,  finite type, exponential decay, Gaussian decay,
OPF operators, polynomial model, weakly pseudoconvex domain}
\subjclass[2000]{32W50 (Primary), 32W30, 32T25}

\maketitle
% main text
%%%%%%%%%%%%%%
%
%	SECTION: INTRODUCTION
%
%%%%%%%%%%%%%%
\section{Introduction}\label{sec:intro}
Let $p:\C\to\R$ be a subharmonic, nonharmonic polynomial
and $\tau\in\R$ a parameter. If $z = x_1 + i x_2$ and
$\frac{\p}{\p\z} = \frac 12 \left(\frac{\p}{\p x_1} 
+ i \frac{\p}{\p x_2}\right)$, define $\Zbstp $
to be the operator 
\[
\Zbstp  = \frac{\p}{\p\z} +  \tau \frac{\p p}{\p\z},
\] 
and let $\Zstp  = -\Zbstp^* = \frac{\p}{\p z} - \tau \frac{\p p}{\p z}$ 
be the negative of the
formal $L^2$-adjoint of $\Zbstp $. If  $ \Boxtp = -\Zbstp\Zstp $, 
then our goal is to  understand the
heat equation:
\begin{equation}\label{eq:he i1}
\begin{cases} {\displaystyle \frac{\p u}{\p s}} + \Boxtp u=0\vspace*{.1in}\\
u(0,z) = f(z).\end{cases}
\end{equation}
We show that the solution $u(s,z)$  of \eqref{eq:he i1}
can be realized as an integral against  a distributional kernel. 
Specifically, we will find a solution to \eqref{eq:he i1}
of the form:
\[
u(s,z) = \int_{\C} \Htp(s,z,w) f(w)\, dw,
\]
and our goal is to understand the regularity and pointwise bounds of $\Htp$ and its
derivatives.
We show that the heat kernel $\Htp(s,z,w)$ is smooth away from the diagonal
$\{(s,z,w) : s=0,\ z=w\}$ and our main result is that:
%
%	THEOREM: POINTWISE BOUNDS ON H(S,Z,W)
%	
\begin{thm}\label{thm:bounds on H}
Let $p$ be a subharmonic, nonharmonic polynomial and $\tau>0$ a parameter.
If $n\geq0$ and $Y^\alpha$ is a product of $|\alpha|$ operators 
$Y = \Zbstp$ or $\Zstp$ when acting in $z$ and $\overline{(\Zstp)}$ or $\overline{(\Zbstp)}$ 
when acting in $w$,
there exist constants $c,c_1>0$ independent of $\tau$ so that
\[
\left|\frac{\p^n}{\p s^n}Y^\alpha  H_{\tau p}(s,z,w)\right|
\leq c_1 \frac{1}{s^{n + \frac 12|\alpha|+1}}
e^{-\frac{|z-w|^2}{32s}} e^{-c \frac{s}{\mu(z,1/\tau)^2}}
e^{-c \frac s{\mu(w,1/\tau)^2}}.
\]
Also, $c$ can be taken with no dependence on $n$ and $\alpha$.
\end{thm}
$\mu(z,\delta)$ is a size function from the Carnot-Carath\'eodory geometry on polynomial models
defined in Section \ref{sec:notation}. As discussed below, in light of the work
of Kurata \cite{Kurata00}, we give new estimates for the long time behavior
of $\Htp$.  The smoothness of $\Htp$ is expected from the work of Nagel and Stein
\cite{NaSt00} and Christ \cite{Christ91}, 
though the estimates for the derivatives of $\Htp$ are new. Moreover, 
as a consequence of \cite{Christ91, NaSt00} and Fu and Straube \cite{FuSt02}, we
expect the results on the heat equation to have applications to partial differential
equations in several complex variables. In fact, we do obtain such applications which
we now describe.
 
The operators $\Zbstp $ and $\Boxtp$ arise in both problems in one complex variable 
and several complex variables.
As detailed below, $\Zbstp$ is a natural operator to consider when studying the weighted $\dbar$-equation
in $\C$ and the $\dbar_b$-problem on polynomial models in $\C^2$. Also, it turns out that the eigenvalues of
$\Box_{np}$ as $n\to\infty$ are important to understand the compactness of 
$\dbar$-Neumann operator on certain classes of Hartogs domains in $\C^2$.

%%%%%%%%%%%%%
%subsection THE DBAR PROBLEM ON WEIGHTED SPACES
%%%%%%%%%
\subsection{$\dbar$ on Weighted $L^q$ Spaces in $\C$}\label{subsec:weighted dbar}

The interest in the weighted $\dbar$-problem in $\C$ begins with H\"ormander's solution
of the inhomogeneous Cauchy-Riemann equations on pseudoconvex domains in $\C^n$
\cite{Hor65}. A crucial estimate in H\"ormander's work is that for  $\Omega\subset\C$ with
$\diam\Omega\leq 1$, there is a solution $u$ to 
$\dbar u=f$ in $L^2(\Omega,e^{-2p})$ satisfying the estimate
$\int_{\Omega} |u|^2 e^{-2p}\, dz \leq  \int_{\Omega} |f|^2 e^{-2p}\, dz$.
Forn{\ae}ss and Sibony \cite{FoSi91} generalize H\"ormander's weighted $L^2$ estimate to
$L^q$, $1<q\leq 2$. They show
$\dbar u=f$ has a solution satisfying $\left(\int_{\Omega} |u|^q e^{-2p}\, dz\right)^{\frac 1q} 
\leq \frac C{p-1}\left(\int_{\Omega} |f|^q e^{-2p}\, dz\right)^{\frac 1q}$.
They also show that the estimate fails if $q>2$. Berndtsson \cite{Ber92} builds on the 
work of Forn{\ae}ss and Sibony by showing an $L^q$-$L^1$ result. He shows that if
$1\leq q < 2$, then $\dbar u=f$ has a solution so that
$\left(\int_{\Omega} (|u|^2 e^{-p})^q\, dz\right)^{\frac 1q} 
\leq  C_p\int_{\Omega} |f| e^{-p}\, dz$.
Berndtsson also proves a weighted $L^\infty$-$L^q$ estimate when $q>2$.

In \cite{Christ91}, Christ recognizes that it is possible to study the $\dbar$-problem in 
$L^2(\C,e^{-2p})$ by working with a related operator in the unweighted space $L^2(\C)$.
If $\dbar\tilde u =\tilde  f$ 
and both $\tilde u = e^{p} u$ and $\tilde  f= e^pf$ are in $L^2(\C,e^{-2 p})$, then
$\frac{\p \tilde u}{\p\z} = \tilde f \Longleftrightarrow e^{-p} \frac{\p}{\p\z} e^p u =f$.
However, $e^{-p} \frac{\p}{\p\z} e^p u = \ZZ_p  u$, so the $\dbar$-problem on 
$L^2(\C,e^{-2p})$ is equivalent to
the $\ZZ_p $-problem, $\ZZ_p  u=f$,
on $L^2(\C)$. 
Christ solves the $\ZZ_p$-equation $\ZZ_p  u  =f$ 
in $L^2(\C)$. 
Christ proves that $G_{p}  = \Box_p^{-1}$ is a well-defined, bounded, linear operator on $L^2(\C)$.
$R_{p}  = Z_p  G_p $ is the relative
fundamental solution of $\ZZ_p $, i.e. the operator $R_p$ satisfies  $\ZZ_p  Rf = (I-S_p)f$
where $S_p$ is the projection of $L^2(\C)$ onto the $\ker \ZZ_p$. 
He shows that $G_{p} $ and $R_{p} $ can be realized
as fractional integral operators with kernels $G_{p} (z,w)$ and $R_{p} (z,w)$, respectively, and
he finds pointwise upper bounds on the kernels $G_p(z,w)$ and $R_p(z,w)$. 

Berndtsson \cite{Ber96} also solves $\ZZ_p u =f$ for $p$ subharmonic,
but Berndtsson solves the problem on $L^2(\Omega)$
where $\Omega\subset\C$ is a
smoothly bounded domain. Like Christ, he expresses
his $L^2$-minimizing 
solution via a fractional integral operator, though unlike Christ, his analysis is derived through
functional analysis and a careful study of Kato's inequality: 
$\triangle |\alpha| \geq \triangle p |\alpha| - 4|\Box_p \alpha|$
where $\alpha \in C^2(\Omega)$. Berndtsson views $\Box_p$ as a
Schr\"odinger operator. Specifically, if $z=x_1+ix_2$, then
\[
2\Box_p = \frac 12 (-i\nabla-a)^2+V
\]
where $a = (-\frac{\p p}{\p x_2},\frac{\p p}{\p x_1})$ and $V = \frac 12 \triangle p$. 
Expressed in this form, $2\Box_p$ is said to be a
Schr\"odinger operator with magnetic potential $a$ and
electric potential $V$. 
We use this  representation of $\Box_p$ in the proof of Theorem \ref{thm:gaussian}.

%%%%%%%%
% SECTION: PSEUDOCONVEX DOMAINS AND DBAR_B
\subsection{Polynomial models and Hartogs Domains in $\C^2$}
Now that we have established the connection between the weighted $\dbar$-equation
in $\C$ and the operators $\ZZ_p $ and $Z_p$, we now turn to the study of
$\dbar_b$-problem on polynomial models in $\C^2$ 
and their connection with the operators $\ZZ_p $ and $Z_p $.
A polynomial model $M_p$ is the boundary of the unbounded weakly pseudoconvex domain 
$\Omega_p = \{(z_1,z_2)\in \C^2: \Imm z_2 > p(z_1)\}$ 
where $p$ is a subharmonic, nonharmonic polynomial.
Observe that the boundary
$M_p\cong\C\times\R$ and the $(0,1)$-form $\dbar_b$ can be identified with the vector field
$\bar L = \frac{\p}{\p \z_1} - 2i \frac{\p p}{\p\z_1}\frac{\p}{\p \z_2}$.
Under the isomorphism, $\dbar_b$ (defined on $M_p$) becomes
the vector field (still called $\bar L$ by an abuse of notation)
\[
\bar L = \frac{\p}{\p\z} - i\frac{\p p}{\p\z}\frac{\p}{\p t}
\]
defined on $\C\times\R$.
There are a number of approaches
that one can take to study the $\bar L$-problem. One is to take a
partial Fourier transform in $t$ because $\bar L$ is translation invariant.
Under the partial Fourier transform, the vector field
$\bar L$ becomes $\Zbstp =  \frac{\p}{\p\z} + \tau\frac{\p p}{\p\z}$, which we regard
as a one-parameter family of
differential operators on $\C$ indexed by $\tau$. 
Thus, questions about the $\dbar_b$-complex on $M$
are closely connected with the $\dbar$-equation on weighted $L^2$-spaces in $\C$.

To analyze operators on Hartogs domains in $\C^n$, 
mathematicians have recognized that it is often enough
to understand weighted operators on the base space and reconstruct the original operator
via Fourier series \cite{Ber94,FoSi91, Li89}. Recently, on a class of Hartogs domains  $\Omega\subset\C^2$,
Fu and Straube \cite{FuSt02, FuSt04} establish an equivalence between the compactness of the
$\dbar$-Neumann operator and the blowup of the minimal eigenvalue of $\Boxtp$ as $\tau\to\infty$.
Christ and Fu \cite{ChFu05} use the work of Fu and Straube  to show the equivalence of:
compactness of the solving operator of the $\dbar$-Neumann Laplacian, 
compactness of the solving operator of Kohn Laplacian $\Box_b$, and
$b\Omega$ satisfying property $(P)$.

% subsection: heat equations!
\subsection{Heat semigroups and heat kernels}\label{subsec:Box heat equations}
Like Christ, we are interested in inverting $\Boxtp =-\Zbstp\Zstp$. 
For an alternative to Christ's approach,  we can look at the heat semigroup $e^{-s\Boxtp}$ 
and integrate in $s$.  Formally, 
$u = e^{-s\Boxtp}[f]$ solves the heat equation \eqref{eq:he i1}
and  inverts $\Boxtp$ since 
\begin{equation}\label{eq:formal box inverse}
\int_0^\infty e^{-s\Boxtp}\, ds = \Boxtp^{-1}
\end{equation}
and $u(0,z) = e^{-0\Boxtp}[f](z) = f(z)$.

On $M$, Nagel and Stein \cite{NaSt00} investigate the heat semigroup $e^{-s\Box_b}$
to solve the heat equation $\frac{\p u}{\p s} + \Box_b u = $ with initial condition $u(0,z)=f(z)$. 
Their goal is to use estimates of the heat semigroup on $M\cong\C\times\R$ 
to understand $\Box_b$ in a product setting \cite{NaSt04}.
Nagel and Stein define $e^{-s\Box_b}$ with 
the spectral theorem use the Riesz Representation Theorem to write $e^{-s\Box_b}[f](\alpha) = \int_{\C\times\R}H(s,\alpha,\beta)f(\beta)\, d\beta$.
$H$ is a distributional kernel with a nonintegrable singularity when $s=0$ and $\alpha=\beta$,
and $H$  is smooth off of the diagonal. They also
obtain pointwise
estimates on $H(s,\alpha,\beta)$ and its derivatives. A fundamental tool in their argument
is the class of nonisotropic smoothing (NIS) operators \cite{NaRoStWa89,NaSt00}. 

A motivation for this work is
to solve the problem of Christ, i.e. invert
$\Boxtp$ and find pointwise estimates on $G_{\tau p}(z,\zeta)$, using the heat
semigroup $e^{-s\Boxtp}$  method and \eqref{eq:formal box inverse}.
In addition, understanding the heat
equation \eqref{eq:he i1} is an interesting question in its own right.
We follow the ideas of \cite{NaSt00} to prove the existence and regularity of $\Htp$. 
Our substitution for their NIS operators are the one-parameter 
families (OPF) of operators defined
in \cite{Rai05f}.
There is an obstruction, however,
to using the techniques of Nagel and Stein in this setting. 
Due to the partial Fourier transform, it 
appears that we cannot scale in the transformed variable. 
Losing the ability to scale in any variable dooms the scaling
argument of Nagel and Stein. We find
other techniques which allow us to bound the heat kernel and its derivatives with
better decay than the scaling argument would have given.

%%%%%%%%%%%%
%
%	SUBSECTION: PUTTING THE RESULT IN CONTEXT AND FUTURE DIRECTIONS
%
%%%%%%%%%%%%%
\subsection{Discussion of Theorem \ref{thm:bounds on H}}\label{sec:future}

The proof of Theorem \ref{thm:bounds on H} has two steps. First,
we show that $e^{-s\Boxtp}$ is an integral operator with
kernel $H_{\tau p}(s,z,w)$ that is smooth away from $\{(s,z,w) : z=w\text{ and }s=0\}$.
To do this, we use the ideas of \cite{NaSt00} to
develop properties of OPF operators defined in
\cite{Rai05f}.
From there, still following  \cite{NaSt00},
we use the spectral theorem and $L^2$-methods to prove smoothness
of $H_{\tau p}(s,z,w)$.

The second step of our analysis is to prove pointwise estimates on $H_{\tau p}(s,z,w)$ and
its derivatives. This is the content of
Theorem \ref{thm:bounds on H}, and the proof has two stages. In the first stage,  
we  write $2\Boxtp$ as a Schr\"odinger operator, similarly to
Berndtsson \cite{Ber96}. We 
use the Feynman-Kac-It\^o formula \cite{Simon79} to show Gaussian
decay for $H_{\tau p}(s,z,w)$. We show the time decay of $\Htp(s,z,w)$ with an 
$L^2$-energy argument.

The goal of the second stage  is to prove pointwise bounds on the derivatives
$\frac{\p^n}{\p s^n}Y^\alpha  H_{\tau p}(s,z,w)$. The idea is to prove a local $L^2$-bound for 
$\frac{\p^n}{\p s^n}Y^\alpha  H_{\tau p}(s,z,w)$ and its derivatives and pass to a local
$L^\infty$-bound using either
a Sobolev embedding-type
result, Theorem \ref{thm:sobolev 1},
or the subsolution estimation from  Kurata \cite{Kurata00}, 
Lemma \ref{lem:kur sub}. The arguments rely heavily on OPF operators and their ability
to commute with derivatives.

Kurata studies heat kernels in $\R^n$ for Schr\"odinger operators
of the form $L = (-i\nabla - a)^2 + V$  where $a\in C^1$ and $V\in L^q_{loc}(\R^n)$, 
$V\geq 0$. His conditions on $a$ and $V$ are more
general than what we consider, and he 
proves continuity of the heat kernel. If $\deg p =2m$, Kurata shows the bound
$|H_{\tau p}(s,z,w)|\leq \frac{C}{s} 
e^{-c_2\frac{|z-w|^2}{s}} e^{-c_3 (\frac{s}{\mu(z,1/\tau)^2})^{1/{2m}}}$, a 
weaker result than ours. The proof of 
Theorem \ref{thm:bounds on H} exploits the specific structure of $\Boxtp$ and
does not seem to generalize to Kurata's more general operators.

By integrating in $s$, the pointwise estimates on $H_{\tau p}(s,z,w)$ allow us to 
recover estimates on the 
fundamental solution of $\Boxtp$
and compare our work to Christ \cite{Christ91}. If $G_{\tau p}(z,w)$ is the  fundamental solution to
$\Boxtp$, we show the decay:
% Corollary: bounds on G_{tau p}^{-1}
\begin{cor}\label{cor:bounds on Box inv} 
Let $G_{\tau p}(z,w)$ be the integral kernel of the fundamental solution
for $\Boxtp^{-1}$. If $X^\alpha$ is a product of $|\alpha|$ operators of the form
$X^j = \Zbstp, \Zstp$ if acting in $z$ and $\overline{(\Zbstp)}, \overline{(\Zstp)}$ if acting in $w$, then
there exists constants $C_{1,|\alpha|}, C_{2}>0$
so that if $\tau>0$,
\[
|X^\alpha G_{\tau p}(z,w)| \leq C_{1,|\alpha|} \begin{cases}
\log\Big(\frac{2\mu(z,1/\tau)}{|z-w|}\Big)   & |z-w| \leq \mu(z,\tfrac 1\tau),\ |\alpha|=0\vspace*{6pt} \\
|z-w|^{-|\alpha|} & |z-w| \leq \mu(z,1/\tau),\ |\alpha|\geq 1\vspace*{6pt} \\
\displaystyle\frac{e^{-C_2 \frac{|z-w|}{\mu(z,1/\tau)}}e^{-C_2 \frac{|z-w|}{\mu(w,1/\tau)}}}
{\mu(z,1/\tau)^{|\alpha|}} 
 & |z-w| \geq \mu(z, \tfrac 1\tau). \end{cases}
\]
Also, $C_2$ does not depend on $\alpha$.
\end{cor}
Near the diagonal, the estimates of Corollary \ref{cor:bounds on Box inv} agree with the bounds of
$X^\alpha G_{p}(z,w)$ computed by Christ. Away from the diagonal, the bounds of Christ are governed by a
metric equivalent to $d\rho^2 = \frac{1}{\mu(\cdot, 1)^2} ds^2$ where $ds$ is the Euclidean metric. Christ
shows that for some $\ep$, $|X^\alpha G_p(z,w)| \les e^{-\ep \rho(z,w)}$. As shown
in Appendix \ref{sec:est equiv}, in the cases that
the author can compute, the two estimates agree. It would be interesting to determine under
what circumstances the two estimates agree or disagree.

In $\R^n$, $n\geq 3$, Shen \cite{Shen99} obtains  estimates for the decay of the fundamental
solution of $-\triangle + V$, $V$ is a nonnegative Radon measure. 
Interestly, his estimates are sharp even though they are higher
dimensional versions of  Christ's estimates which are not
sharp. This signifies there is additional structure in the 
special relationship between $a$ and $V$ in the magnetic
Schr\"odinger operator $\Boxtp$. 

Once we have estimates for all $\tau\in\R$,
estimates on $\Htp(s,z,w)$ will have many applications to questions in several
complex variables. 
Pointwise
estimates for $\Htp(s,z,w)$ and its derivatives when $\tau <0$  
is the subject of \cite{Rai06}.
A difficulty lies in the fact that techniques from parabolic operator theory and 
quantum mechanics do not seem to work.  
In the Schr\"odinger operator representation,
$2\Boxtp = \frac 12 (-i\nabla -  a)^2 +  V$, $V\leq 0$ and unbounded. Writing
$\Boxtp$ as a parabolic operator, this means the unbounded $0^{\text{th}}$ order term
may not be positive.
We plan to use our estimates of $\Htp(s,z,w)$
to prove exponential decay for the heat kernel of \cite{NaSt00}, an improvement
over the rapid decay shown by Nagel and Stein. 
We also hope to use the OPF operator and heat kernel results to 
build on the work of \cite{NaSt03} by  proving pointwise estimates  on the heat
kernel on the boundary of decoupled domains in $\C^n$, i.e. domains of the form 
$M  = \{(z_1, \ldots, z_n) : \Imm z_n  < \sum_{j=1}^{n-1} p_j(z_j)\}$ where $p_j$ are nonharmonic,
subharmonic polynomials.

\textbf{Acknowledgements.} 
I would like to thank Alexander Nagel for his support and guidance during this
project.

%%%%%%%%%%%%%
%
%
%	SECTION: DEFINITIONS AND RESULTS
%
%%%%%%%%%%%%%
\section{Notation and Definitions}\label{sec:notation}

% SUBSECTION: NOTATION
\subsection{Notation For Operators on $\C$}\label{subsec:notationdown}
For the remainder of the paper,  let $p$ be a  subharmonic, nonharmonic 
polynomial.
It will be important for us to write $p$ centered around an arbitrary point $z\in\C$, and 
we set:
\begin{equation}\label{eq:Ajk}
\A{jk}z = \frac{1}{j! k!}\frac{\p^{j+k}p}{\p z^j\p\z^k}(z).
\end{equation} 
We need the following functions two ``size" functions to write down 
the size and cancellation conditions 
for both OPF operators and NIS operators. Let
\begin{align}\label{eq:Lam}
\Lambda(z,\delta) &= \sum_{j,k\geq 1}  |\A{jk}z |\delta^{j+k}
\intertext{and}
\label{eq:mu}
\mu(z,\delta) &= \inf_{j,k\geq 1} \left|\frac{\delta}{\A{jk}z}\right|^{1/(j+k)}.
\end{align}
$\Lambda(z,\delta)$ and $\mu(z,\delta)$ are geometric objects from the Carnot-Carath\'eodory
geometry  developed by Nagel, Stein, and Wainger 
\cite{Na86, NaStWa85}. The functions also arise
in the analysis of  magnetic Schr\"odinger  operators with electric potentials \cite{Kurata00, Shen96,
Shen99}. 
It follows 
$\mu(z,\delta)$ is an approximate inverse to $\Lambda(z,\delta)$. This means that if $\delta>0$,
\[
\mu\big(z,  \Lambda(z,\delta)\big) \sim \delta
\text{ and }
\Lambda\big(z,\mu(z,\delta)\big) \sim \delta.
\]
We use the notation $a\les b$ if $a\leq C b$ where $C$ is a constant that may depend on the
dimension 2 and the degree of $p$. We say that $a\sim b$ if $a\les b$ and $b\les a$.
 
Denote the ``twist" at $w$, centered at $z$ by 
\begin{equation}
T(w,z) = - 2\Imm \left( \sum_{j \geq 1} \frac{1}{j!} \frac{\p^j p}{\p z^j}(z) (w-z)^j\right).
\label{eq:twist}
\end{equation}

Also associated to a polynomial $p$ and the parameter $\tau\in\R$ are the weighted differential
operators 
\begin{align*}
\Zbstpz &= \frac{\p}{\p \z} + \tau \frac{\p p}{\p \z}
= e^{-\tau p}\frac{\p p}{\p \z}e^{\tau p},
& \Zstpz &= \frac{\p}{\p z} - \tau \frac{\p p}{\p z}  
= e^{\tau p}\frac{\p p}{\p z}e^{-\tau p}\\
\intertext{and}
\Wbstpw &= \frac{\p}{\p \w} -  \tau \frac{\p p}{\p \w}
= e^{\tau p}\frac{\p p}{\p \w}e^{-\tau p},
& \Wstpw &= \frac{\p}{\p w} +  \tau \frac{\p p}{\p w}
= e^{-\tau p}\frac{\p p}{\p w}e^{\tau p}.
\end{align*}
$\Zbstpz$ and $\Zstpz$ arise naturally as described above. The need for 
$\Wbstpw$ and $\Wstpw$ is explained below.
 
We think of $\tau$ as fixed and the operators $\Zbstpz$, 
$\Zstpz$, $\Wbstpw$, and $\Wstpw$ as acting
on functions defined on $\C$.
Also, we will omit the variables $z$ and $w$ from 
subscripts when the application is unambiguous. Observe that
$\overline{(\Zstp)} =  \Wbstp$ and $\overline{(\Zbstp)} = \Wstp$.
We let $X_1$ and $X_2$ denote the ``real" and ``imaginary" parts of $Z$, that is,
\begin{align*}
X_1 &= \Zstp + \Zbstp = \frac{\p}{\p x_1} + i \tau \frac{\p p}{\p x_2}, &
X_2 &= i(\Zstp - \Zbstp) = \frac{\p }{\p x_2} - i \tau \frac{\p p}{\p x_1}.
\end{align*}
Analogously to $X_1$ and $X_2$, define 
\begin{align*}
U_1 &= \Wstp + \Wbstp = \frac{\p}{\p x_1} - i \tau \frac{\p p}{\p x_2}, &
U_2 &= i(\Wstp - \Wbstp) = \frac{\p }{\p x_2} +i \tau \frac{\p p}{\p x_1}.
\end{align*}

We need to establish notation for adjoints. 
If $T$ is an operator (either 
bounded or closed and densely defined) on a Hilbert space with inner product
$\big(\,\cdot\, ,\cdot\,\big)$, let $T^*$ be the Hilbert space adjoint of $T$. This means 
that if
$f\in\Dom T$ and $g\in \Dom{T^*}$, then $\big(Tf,g\big) = \big(f,T^*g\big)$. 
If $U$ is an unbounded domain in some Euclidean space,
$T$ is an operator acting on $C^\infty_c(U)$ or $\mathcal{S}(U) = \{\vp\in C^\infty(U) : 
\vp \text{ has rapid decay}\}$, then we denote $T\sh$ as the adjoint in the sense of distributions.
This means if $K$ is a distribution or a Schwartz distribution, then
$\langle T\sh K, \vp\rangle = \langle K, T \vp \rangle$. Note that if $T$ is not 
$\R$-valued, $T^* \neq T\sh$.
It follows easily that 
\[
\Zbstp\sh = -\Wbstp \qquad\text{and}\qquad \Zstp\sh = -\Wstp.
\]
Finally, let $\Ms = \T$.

% SUBSECTION: DEFINITION OF THE FAMILY OF OPERATORS
\subsection{Definition of OPF Operators}\label{subsec:OPF operator}
\renewcommand{\labelenumi}{(\alph{enumi})}
We use the definition from \cite{Rai05f}.
\renewcommand{\labelenumi}{(\alph{enumi})}
We say that  $\Ts$ is a \emph{one-parameter family (OPF)
of operators}  of order $m$ with respect to the polynomial $p$ if the 
following conditions hold:
\begin{enumerate}  

\item There is a  function 
$\Ks\in C^\infty\Big(\big((\C\times\C) \setminus\{z=w\}\big) \times (\R\setminus\{0\})\Big)$ 
so that for fixed $\tau$,  
$\Ks$ is a distributional kernel, i.e. if $\vp,\psi \in C^\infty_c(\C)$ and 
$\supp\vp \cap \supp \psi = \emptyset$, then
$\Ts[\vp]\in (C^\infty_c)'(\C)$ and  
\[
 \langle \Ts[\vp](\cdot),\psi\rangle_{\C} = \iint_{\C\times\C} \Ks(z,w)\vp(w) \psi(z)\, dw dz.
\]

\item \label{it:approx}
There exists a family of functions $\Kse(z,w)\in C^{\infty}(\C\times\C\times\R)$ so
that if $\vp\in C^\infty_c(\C\times\R)$,
\[
 \Kse[\vp]_{\C\times\R}(z,\tau) = \int_{\C\times\R}\vp(w,\tau)\Kse(z,w)\,dw d\tau
\]
and $\lim_{\epsilon\to0} \Kse[\vp]_{\C\times\R}(z) = \Ks[\vp]_{\C\times\R}(z)$ in 
$(C^\infty_c)'(\C\times\R)$.  \\
All of the additional conditions are assumed to apply to the kernels 
$\Kse(z,w)$ uniformly in $\ep$.

\item 
Size Estimates. If $\Ystp^{J}$ is a product of $|J|$ operators of the form
$\Ystp^j =\Zstpz$, $\Zbstpz$,  $\Wstpw$, $\Wbstpw$, or
$\Ms$ where $|J| = \ell +n$ and
$n = \#\{j : \Ystp^j = \Ms\}$, for any $k\geq 0$ there exists a constant $C_{\ell,n,k}$ so that
{\small
\begin{equation}\label{it:size}
\left| \Ystp^{J}  \Kse(z,w)\right| \leq \frac{C_{\ell,n,k}}{|\tau|^n} 
\frac{|z-w|^{m-2-\ell}}{|\tau|^{k}\Lambda(z,|w-z|)^k} \quad \text{if}\quad
\begin{cases} &m<2 \\ &m=2,\ k\geq 1 \\ &m=2, |w-z| > \mu(z,\tfrac 1\tau) \end{cases} 
\end{equation}
}
Also, if $m=2$ and $|w-z| \leq \mu(z,\tfrac 1\tau)$, then
\begin{equation}\label{it:size2}
\left| \Ms^n \Kse(z,w) \right| \leq C_n 
\begin{cases} \log\left( \frac{2\mu(z,\tfrac 1\tau)}{|w-z|}\right)  &n=0 \\
|\tau|^{-n} &n\geq 1 \end{cases}
\end{equation}

\item 
Cancellation in $w$.  If $\Ystp^{J}$ is a product of $|J|$ operators of the form
$\Ystp^j = \Zstpz$, $\Zbstpz$, $\Wstpw$, $\Wbstpw$, or $\Ms$ where $|J| = \ell +n$ and
$n = \#\{j : \Ystp^j = \Ms\}$, for any $k\geq 0$ 
there exists a constant $C_{\ell,n,k}$ and $N_\ell$ 
so that for $\vp\in C^\infty_c (D(z_0,\delta))$, {\small
\begin{align}\label{it:cancel w}
& \sup_{z\in\C} \left| \int_{\C}  \Ystp^{J}\Kse(z,w)\vp(w)\,dw\right|\nn \\
& \leq \frac{C_{\ell,n,k}}{|\tau|^n} 
 \begin{cases} 
   {\displaystyle\delta^2
  \Big(\log\big(\tfrac{2\mu(z,1/\tau)}{\delta}\big) \|\vp\|_{L^\infty} \hspace{-2pt}
  + \hspace{-8.5pt}\sum_{1 \leq |I|\leq N_0} \delta^{|I|}\|\Xstp^I\vp(w)\|_{L^\infty}\Big)} 
  &\atopp{\displaystyle \delta < \mu(z,\tfrac 1\tau) 
\text{ and}}{\displaystyle m=2,\ell=0} \\
 {\displaystyle \frac{\delta^{m-\ell}}{|\tau|^{k}\Lambda(z,\delta)^k} \sum_{|I| \leq N_\ell}
 \delta^{|I|} \left\| \Xstp^I \vp\right\|_{L^\infty(\C)}} & \text{otherwise}
\end{cases}
\end{align}
}
 where $\Xstp^I$ is composed solely of  $\Zstp$ and $\Zbstp$.

\item Cancellation in $\tau$. If $\Xstp^{J}$ is a product of $|J|$ operators of the form
$\Xstp^j = \Zstpz,\ \Zbstpz$ or $\Wstpw,\ \Wbstpw$ 
and $|J| = n$,
there exists a constant $C_{n}$ so that
\begin{equation} \label{it:cancel tau}
 \int_\R\Xstp^{J}\left( e^{ i \tau t}\Kse(z,w)\right) \, d\tau \leq 
 C_n \frac{\mu(z,t+T(w,z))^{m-n}}{\mu(z,t+T(w,z))^2 |t+T(w,z)|}.
\end{equation}

\item \label{it:adjoint} Adjoint.
Properties (a)-(e) 
also hold for the adjoint operator $\Ts^*$ whose distribution kernel is given
by $\overline{\Kse(w,z)}$
\end{enumerate}

The following results from \cite{Rai05f} are essential tools in the proof of Theorem
\ref{thm:bounds on H}.
%
%	THEOREM: L^P BOUND
%
\begin{thm}\label{thm:Lp bound}
If $\Ts$ is an OPF operator of order 0, then $\Ts$, $\Ts^*$ are bounded operators from
$L^q(\C)$ to $L^q(\C)$, $1<q<\infty$, with a constant independent of $\tau$ but
depending on $q$.
\end{thm}
%
%	THEOREM: ONE TO ONE CORRESPONDENCE
%
\begin{thm}\label{thm:correspondence}
Given a subharmonic, nonharmonic polynomial $p:\C\to\R$, 
there is a one-to-one correspondence between OPF operators
of order $m\leq 2$ with respect to $p$
and NIS operators of order $m\leq 2$ on the polynomial model 
$M_p = \{(z_1,z_2)\in\C^2 : \Imm z_2 = p(z_1)\}$. The correspondence is given by a partial
Fourier transform in $\Rre z_2$.
\end{thm}

There are multiple definitions of NIS operators (e.g. \cite{NaRoStWa89, NaSt00}).
This equivalence is with the definition in \cite{NaRoStWa89}.

\section{The Heat Equation and Smoothness of the Heat Kernel}\label{sec:heat smoothness}

For the remainder of the work, we will primarily  be concerned with
inverting the ``Laplace" operator 
\[
\Boxtp = -\Zbstp \Zstp.
\]
via the heat semigroup $e^{-s\Boxtp}$. 
We assume that $\tau>0$ and  define the heat operator
\[
\H = \frac{\p}{\p s} + \Boxtp.
\]
Given a function $f$ defined on $\C$, we study the initial value
problem of finding smooth $u :(0,\infty)\times \C \to \C$ so that
\begin{equation}\label{eq:he1}
\begin{cases} \H[u](s,z) =0 &s>0,\ z\in\C\\
\displaystyle \lim_{s\to 0} u(s,\cdot) = f(\cdot) &\text{with convergence in an appropriate norm.}
\end{cases}
\end{equation}

%We will show that the solutions $u$ 
%defined on 
%$(0,\infty) \times \C$ is given by semigroup of operators
%\[
%u(s,z) = e^{-s \Boxtp}[f](z)
%\]
%We also prove  the existence of a function $\Htp(s,z,w)$ which is smooth
%away from the diagonal $\{(s,z,w) : s=0 \text{ and } z=w\}$ and has the
%property  that the heat semigroup
%$e^{-s \Boxtp}$ can be written
%\[
% e^{-s \Boxtp}[f](z) = \int_{\C} f(w) \Htp(s,z,w)\, dw.
%\]

Let $\alpha$ be a multiindex. We let $X^\alpha$ be a product of $|\alpha|$ operators of the form
$X = X_1$ or $X_2$. Similarly, $U^\alpha$ is a product of $|\alpha|$ operators of the form
$U = U_1$ or $U_2$.

%%%%%%%%%%%%%%%%%%%%%%%
%
%		SECTION: SPECTRAL THEOREM
% 
%%%%%%%%%%%%%%%%%%%%%%%
\renewcommand{\labelenumi}{(\alph{enumi})}
\section{The heat semigroup $e^{-s\Boxtp}$ on $L^2(\C)$} \label{sec:semigroups}
We know that $\Zbstp$ and $\Zstp$ are closed, densely defined operators on $L^2(\C)$. 
As in Nagel-Stein \cite{NaSt00}, the spectral theorem for unbounded operators (see \cite{Rudin91})
proves:
% Theorem: spectral theorem/semigroups
\begin{thm} \label{thm:semigroups} $\Boxtp$ is the
infinitesimal generator of $e^{-s\Boxtp}$, a strongly continuous semigroup of bounded operators on 
$L^2(\C)$  for $s>0$. For $f\in L^2(\C)$, the following hold:
\begin{enumerate}
 \item \label{item:lim e^-sbox = f}
  $\displaystyle \lim_{s\to 0}\|e^{-s\Boxtp} f - f \|_{L^2(\C)} =0;$
 
 \item \label{item:e^-sbox contraction}
  For $s>0$, these operators are contractions, that is,
\[
 \| e^{-s\Boxtp} f \|_{L^2(\C)} \leq \| f\|_{L^2(\C)};
\]

\item \label{item:e^-sbox -f < f}
  For $f\in\Dom(\Boxtp)$,
\[ 
 \|e^{-s\Boxtp} f - f \|_{L^2(\C)} \leq s \|\Boxtp f\|_{L^2(\C)};
\]

\item \label{item:box^j f bound}
For $s>0$ and all $j$, $\Ran(e^{-s\Boxtp}) \subset \Dom(\Boxtp^j)$. 
Also, $\Boxtp^j e^{-s\Boxtp}$ is a bounded operator on $L^2(\C)$ with
\[
\| \Boxtp^j e^{-s\Boxtp} f \|_{L^2(\C)} \leq \left( \frac je \right) s^{-j} \| f\|_{L^2(\C)}; 
\]

\item\label{item:heat eqn} For any $f\in L^2(\C)$ and $s>0$, the Hilbert space valued function
$u(s) = e^{-s\Boxtp} f$ satisfies
\[
\left(\frac{\p}{\p s} + \Boxtp \right) u(s) =0.
\]
\end{enumerate}
\end{thm}

%%%%%%%%%%%%%%%%%%%
%
%		SECTION: REGULARITY OF THE HEAT KERNELS
%
%%%%%%%%%%%%%%%%%%%
\section{Regularity of the Heat Kernel}\label{sec:regularity}
\begin{def}\label{def:Hs} For each $s>0$, define the bounded operator $\Htp^s:
L^2(\C)\to L^2(\C)$  by
\[
\Htp^{s}[f] = e^{-s\Boxtp}[f].
\]
\end{def}
%
%%%% THEOREM: PROPERTIES OF G(S,Z,W), H(S,Z,W)
%
Our main result on the existence and regularity of $\Htp(s,z,w)$ is the following theorem. 
\begin{thm}\label{thm:H properties}
Fix $\tau>0$. There is a function 
$\Htp \in C^\infty\big((0,\infty)\times\C\times\C\big)$ so that
for all $f\in L^2(\C)$,
\begin{equation}
\Htp^{s}[f](z) = \int_\C \Htp(s,z,w) f(w)\, dw \label{eq:H def}
\end{equation}
Moreover, for each fixed $s>0$ and $z\in \C$, the function $w\mapsto \Htp(s,z,w)$
is in $L^2(\C)$, so the integral defined in equation \eqref{eq:H def}
converges absolutely.  Also,
\begin{enumerate}

\item\label{item:conjugate} $\displaystyle \Htp(s,z,w) = \overline{\Htp(s,w,z)}$;

\item\label{item:heat equ on kernels} For $(s,z,w)\in (0,\infty)\times\C\times\C$,
\[
\left(\frac{\p}{\p s} + \Boxtpz\right)[\Htp](s,z,w) = \left(\frac{\p}{\p s} + \Boxtpw\sh\right)[\Htp](s,z,w) =0;
\]

\item\label{item:Boxz = Boxw} For any integers $j,k\geq 0$,
\[
\Boxtpz^j(\Boxtpw\sh)^k \Htp(s,z,w) = \Boxtpz^{j+k}\Htp(s,z,w) = (\Boxtpw\sh)^{j+k} \Htp(s,z,w);
\]

\item\label{item:L^2 derivatives of H}
For all integers $j$ and multiindices $\alpha, \beta$, the functions
\[
w \longmapsto \frac{\p^j}{\p s^j} X_z^\alpha U_w^\beta \Htp(s,z,w)
\]
are in $L^2(\C)$ and there is a constant $c_{j,\alpha,\beta}$ so that for $R < R_{\tau p}(z)$,
\[
\big\| \frac{\p^j}{\p s^j} X_z^\alpha U_w^\beta \Htp(s,z,\cdot) \big\|_{L^2(\C)} 
\leq  \frac {C_{\alpha,\beta, j}}R s^{-\frac{\alpha+\beta}2-j}(1+ s^{-1});
 \]

\item\label{item:z,w interchange} The conclusions of    (d) 
hold with the roles of $z$ and $w$ interchanged.
\end{enumerate}
\end{thm}

%%%%%%%%%%%%
%
%	Section: FAMILY PROPERTIES
%
%%%%%%%%%%%%%
\section{Properties of OPF Operators}
\label{sec:properties}
To prove Theorem \ref{thm:H properties}, we need to establish
properties of OPF operators. We follow the line of argument for
NIS operators in \cite{NaRoStWa89,NaSt00}.
Since we are working with a fixed polynomial $p$, we omit $\tau p$
from subscripts when the application is clear.

% Lemma: "commuting" X around order 0 operators
\begin{lem}\label{lem:X commuting}
Let $A_\tau$ and $B_\tau$  be order $0$ OPF operators, and let $X = X_1$ or $X_2$.
There exist order 0 OPF operators 
$A_1$, $A_2$, $B_1$ and $B_2$ so that
\begin{align*}
X A_\tau &= A_1 X_1 + A_2 X_2; \\
B_\tau X&= X_1B_1 + X_2 B_2.
\end{align*}
\end{lem}

\begin{proof} We know from results for NIS operators and Theorem \ref{thm:correspondence}
that $X_1^2 + X_2^2$ is invertible
with an inverse $K_\tau$ that is an OPF operator smoothing of order 2. Thus, we can write
\[
X A_\tau =  \Big( X A_\tau K_\tau X_1\Big) X_1 + 
\Big( X A_\tau K_\tau X_2\Big) X_2 = A_1 X_1 + A_2 X_2.
\]
A similar argument proves the result for $B_{\tau}$.
%Similarly,
%\[
%B_\tau X = X_1\Big(X_1 K_\tau B_\tau X \Big) +  X_2\Big(X_1 K_\tau B_\tau X \Big)
%= X_1B_1 + X_1B_2.
%\]
\end{proof}

% Corollary:  "commuting X^\alpha around order 0 operators
\begin{cor}\label{cor:X^alpha commuting} 
Let $A_\tau$ and $B_\tau$  be order $0$ OPF operators and $\alpha$ a multiindex where
$|\alpha| = k\geq 1$. There exist finite sets $I$ and $J$ of multiindices $\alpha_i$, $|\alpha_i| =k$,
and $\beta_j$, $|\beta_j|=k$, respectively so that 
\begin{align*}
X^\alpha A_\tau &= \sum_{\alpha_i\in I} A_i  X^{\alpha_i},\\
B_\tau X^\alpha&= \sum_{\alpha_j\in J} X^{\beta}B_j
\end{align*}
for some order 0 OPF operators $A_i$ and $B_j$.
\end{cor}

\begin{proof} Induction.
\end{proof}
Let $[n]$ denote the greatest integer less than or equal to $n$. The proofs of the following
two propositions are in \cite{NaSt00}.

% Proposition:  changing X's into A \Box.
\begin{prop}\label{prop:X to Box A}
Let $\alpha$ be a multiindex.
\begin{enumerate} 
\item If $|\alpha|$ is even, there exists an order 0 OPF operator
$A_\tau$ so that
\[
X^\alpha = \Boxtp^{\frac{|\alpha|}2} A_\tau.
\]
\item  If $|\alpha|$ is odd, there exist order 0 OPF operators
$A_1$ and $A_2$ 
so that 
\[
X^{\alpha} = \Boxtp^{[\frac{|\alpha|}2]} (X_1 A_1 + X_2 A_2).
\]
\end{enumerate}
\end{prop}

% Proposition:  changing X's into  \Box B.
\begin{prop}\label{prop:X to B Box}
Let $\alpha$ be a multiindex.
\begin{enumerate} 
\item If $|\alpha|$ is even, there exists an order 0 OPF operator 
$B_\tau$  so that
\[
X^\alpha = B_\tau\Boxtp^{\frac{|\alpha|}2};
\]
\item  If $|\alpha|$ is odd, there exist order 0 OPF operators
$B_1$ and $B_2$ 
so that 
\[
X^{\alpha} = \Boxtp^{[\frac{|\alpha|}2]} (B_1X_1 + B_2X_2);
\]
\item Alternatively, if $|\alpha|$ is odd and $X^\alpha = X^\beta X$ where
$X = X_1$ or $X_2$, then there exists an order 0 
OPF operator $B_\tau$ so that
\[
X^\alpha = B_\tau \Boxtp^{\frac{|\beta|}2}X.
\]
\end{enumerate}
\end{prop}

%\begin{proof} The proof is almost identical to the proof of Proposition \ref{prop:X to Box A}.
%\end{proof}

% PROPOSITION: L^2 X BOUND IN TERMS OF THE FUNCTION AND \BOX
\begin{prop} \label{prop: L^2 X bound}
Let $X = X_1$ or $X_2$.  There is a constant $C$ so that if $\vp\in C^\infty_c(\C)$, then
for all $r>0$
\[
\| X[\vp]\|_{L^2(\C)} \leq C ( r\|\Boxtp\vp\|_{L^2(\C)} + r^{-1} \|\vp\|_{L^2(\C)}).
\]
\end{prop}

\begin{proof} 
First, note that $X^* = -X$. Using Proposition \ref{prop:X to B Box}, we compute
\begin{align*}
\| X[\vp]\|^2_{L^2(\C)} = \big( X[\vp], X[\vp] \big)& = -\big(  X^2[\vp],\vp\big)
\leq \left|\big( A_\tau\Boxtp[\vp],\vp\big)\right|  \\
&\leq C \big(r^2\|\Boxtp\vp\|_{L^2(\C)}^2 + r^{-2} \|\vp\|_{L^2(\C)}^2 \big).
\end{align*}
\end{proof}

% COROLLARY: L^2 X^ALPHA BOUNDS IN TERMS OF \VP AND \BOX 
\begin{cor}\label{cor: L^2 X^alpha bound}
Let $\alpha$ be a multiindex. There exists a constant  $C_{|\alpha|}$ so that if
$\vp\in C^\infty_c(\C)$, then 
\[
\| X^\alpha [\vp] \|_{L^2(\C)} \leq C_\alpha 
\sum_{j=0}^{[\frac {|\alpha|}2]+1} \| \Boxtp^j[\vp] \|_{L^2(\C)}.
\]
\end{cor}

\begin{proof} Proof by induction. The base case is Proposition \ref{prop: L^2 X bound}
and the inductive step is a repetition of the argument in the proof of
Proposition \ref{prop: L^2 X bound}.
\end{proof}

% THEOREM: SOBOLEV TYPE 1
We now prove the Sobolev type theorem.
\begin{thm} \label{thm:sobolev 1}
Let 
\[
R_{\tau p}(z) =\inf_{j,k\geq 0} \frac{1}{|\tau\A{jk}{z}|^{\frac 1{j+k}}}.
\]
There is a constant $C>0$ 
so that if $f\in C^\infty (\C)$, $z\in \C$ and $0<R < R_{\tau p}(z)$,
\[
\sup_{D(z,R)}|f| \leq \frac {C}{R} \sum_{|\alpha|\leq 2} R^{|\alpha|} 
\| X^\alpha f \|_{L^2(D(z,2R))}.
\]
Also, if $f\in C^\infty(\C)\cap L^2(\C)$, then
\[
\sup_{D(z,R)}|f| \leq \frac {C}{R}\left( \| f \|_{L^2(D(z,2R))}+ R^{2}\|\Boxtp f \|_{L^2(D(z,2R))}\right).
\]
\end{thm}

\begin{proof}
Let $f\in C^\infty(\C)$ and $z\in\C$. An application of Plancherel's Theorem shows
\begin{equation}\label{eqn:sobo}
\sup_{D(z_0,R)}|f(z)| \leq \frac CR \sum_{|\alpha|\leq 2} R^{|\alpha|} 
\|D^\alpha f(z)\|_{L^2(D(z_0,2R))}.
\end{equation}

%Let $\chi\in C^\infty_c(\C)$, 
%$\chi \equiv 1$ on $D(0,1)$, $0\leq \chi\leq 1$, and $\chi(z) \equiv 0$  if $|z| \geq 2$.
%Let $g(z) = f(z)\chi(z-z_0)$. Then 
%\begin{align*}
%\sup_{D(z_0,1)}|f(z)|  &
%\leq \sup_{\C} |g| \leq \int_{\C} |\hat g(\xi)|\, d\xi
%\leq \|(1+|\xi|^4)^{\frac12} \hat g(\xi)\|_{L^2(\C)} \|(1+|\xi|^4)^{-\frac12}\|_{L^2(\C)}\\
%&\leq C\big(\|\hat g\|_{L^2(\C)} + \||\xi|^2\hat g\|_{L^2(\C)}\big)
%\leq C\big(\|g\|_{L^2(\C)} + \|\triangle g\|_{L^2(\C)}\big) 
%\leq \sum_{|\alpha|\leq 2} \|D^\alpha f\|_{L^2(D(z_0,2))}.
%\end{align*}
%A  change of variables argument shows
%\begin{equation}\label{eqn:sobo}
%\sup_{D(z_0,R)}|f(z)| \leq \frac CR \sum_{|\alpha|\leq 2} R^{|\alpha|} 
%\|D^\alpha f(z)\|_{L^2(D(z_0,2R))}.
%\end{equation}

To pass from ordinary derivatives to products of $\Zstp$ and $\Zbstp$, first observe that if 
$|w-z|<R<R_{\tau p}(z)$,
then
\begin{align*}
\left|\tau\frac{\p^{j+k} p}{\p z^j\p\z^k}(w)\right| 
&= \left|\tau \sum_{\atopp{j'\geq j}{k'\geq k}} 
\frac{1}{(j-j')!(k-k')!} \frac{\p^{j'+k'}p}{\p z^{j'}\p\z^{k'}}(z)
(w-z)^{j'-j}\overline{(w-z)}^{k'-k}\right|\\
&\leq \frac{C \tau}{R^{j+k}}\sum_{\atopp{j'\geq j}{k'\geq k}} 
 \left|\frac{\p^{j'+k'}p}{\p z^{j'}\p\z^{k'}}(z)\right| R^{j'+k'} \leq \frac{C}{R^{j+k}},
\end{align*}
where $C$ does not depend on $p$ or $R$. 
The proof of the first part of the theorem now follows easily since $\frac{\p f}{\p\z}(w)
= \Zbstp f(w) -   \tau \frac{\p p}{\p\z}(w)f(w)$. This means

\[
\left| \frac{\p f}{\p \z}(w)\right| \leq |\Zbstp f(w)| + \frac CR|f(w)|
\]
and similarly for  $\left| \frac{\p f}{\p z}(w)\right|$. 
The estimates for the second derivatives proceed in the same fashion. 
For example, {\small
\begin{align*}
&\left| \frac{\p^2  f}{\p z\p\z} (w)\right|\\
&= \left| \Zbstp\Zstp[f](w) -  \tau\tfrac{\p p(w)}{\p\z}\tfrac{\p f(w)}{\p z}
+\tau^2\tfrac{\p p(w)}{\p\z}\tfrac{\p p(w)}{\p z} f(w)+  \tau\tfrac{\p p(w)}{\p z}\tfrac{\p f(w)}{\p \z}
+ \tau \tfrac{\p^2  p(w)}{\p z\p\z} f(w) \right|\\
&\leq |\Zbstp\Zstp[f](w) | + \frac {C}{R}|\nabla f(w)| + \frac{C}{R^2}|f|.
\end{align*} }
The other second derivatives of $f$ are handled similarly. 
Thus, every term in \eqref{eqn:sobo} is well controlled by
$\Zstp$ and $\Zbstp$ derivatives. 
The proof of the latter part of the theorem follows from Proposition
\ref{prop:X to B Box} and Proposition  \ref{prop: L^2 X bound}.
\end{proof}

\begin{rem} In Theorem \ref{thm:sobolev 1}, 
if $\frac{\p^j p}{\p z^j}(z) = \frac{\p^j p}{\p\z^j}(z) =0$ for
$j= 0,1,\ldots,\deg(p)$, then $R_{\tau p}(z) = \mu(z, 1/\tau)$, 
a fact which will be useful later.
\end{rem}

%%%%%%%%%%%%%%
%
%	PROOF OF THE SMOOTHNESS THEOREM
%
%%%%%%%%%%%%%%
\section{Proof of Theorem \ref{thm:H properties}} \label{sec:smoothnessproof}

To prove Theorem \ref{thm:H properties}, we need some a priori estimates.

% Lemma: X^\alpha G_s X^\beta \vp L^2 estimates

\begin{lem}\label{lem:X^alpha H^{s} X^beta phi L^2 est}
There are constants $C_{\alpha,\beta}$ so that for any multiindices $\alpha$ and $\beta$,
any $s>0$, and $\vp\in C^\infty_c(\C)$,
\[
\| X^\alpha \Htp^{s} [X^\beta \vp] \|_{L^2(\C)} 
\leq C_{\alpha,\beta} s^{-\frac{|\alpha|+|\beta|}2} \|\vp\|_{L^2(\C)}.
\]
\end{lem}

\begin{proof} We first assume that 
$|\alpha|$ and 
$|\beta|$ are even. From Proposition \ref{prop:X to Box A}, there exists 
an order 0 OPF operator $A_\tau$ so that
\[
 \Htp^{s}[X^\beta \vp] = \Htp^{s} \Boxtp^{\frac{|\beta|}2} A_\tau \vp.
\]
Hence, we have by Proposition \ref{prop:X to B Box} and Theorem \ref{thm:semigroups} (d) an
order zero family $B_\tau$ so that
\begin{align*}
\| X^\alpha \Htp^{s} [X^\beta \vp] \|_{L^2(\C)}
& = \| X^\alpha \Boxtp^{\frac{|\beta|}2}\Htp^{s} [A_\tau\vp] \|_{L^2(\C)}\\
& = \| B_\tau \Boxtp^{\frac{|\alpha|}2} \Boxtp^{\frac{|\beta|}2}\Htp^{s} 
[ A_\tau\vp] \|_{L^2(\C)} \\
&\leq C_{\alpha,\beta} s^{- \frac{|\alpha|+|\beta|}2}\|  \vp\|_{L^2(\C)}.
\end{align*}
The $|\alpha|$ and $|\beta|$ odd cases follow easily from the even case, an application
of Proposition \ref{prop:X to Box A} and Proposition \ref{prop:X to B Box}, 
and the following two arguments.
One, if $X$ is  either $X_1$ 
or $X_2$ then from Proposition \ref{prop: L^2 X bound} with $r = s^{\frac12}$,
\[
\| X \Htp^s \vp \|_{L^2(\C)}  
\leq C\Big( s^{\frac 12} \|\Boxtp \Htp^s \vp\|_{L^2(\C)} + s^{-\frac 12} \| \Htp^s \vp\|_{L^2(\C)}\Big )
\leq C s^{-\frac 12}  \|  \vp\|_{L^2(\C)}.
\]
Two, since $X^* = -X$, applying the previous inequality to $\Htp^{s} X\vp$, we have
\begin{align*}
\| \Htp^{s} X\vp \|_{L^2(\C)}^2 &= \Big(\Htp^{s} X\vp, \Htp^{s} X\vp \Big)
= - \Big( \vp, X \Htp^{s}\Htp^{s} X\vp \Big) \\
&\leq \|\vp\|_{L^2(\C)} \|X \Htp^{s} \Htp^{s} X\vp\|_{L^2(\C)}
\leq C s^{-\frac 12} \|\vp\|_{L^2(\C)} \|  \Htp^{s} X\vp \|_{L^2(\C)}.
\end{align*}
\end{proof}

%Lemma: X^alpha H^{s}[\vp] pointwise estimates
\begin{lem}\label{lem:X H_[s] pointwise est}
For $s>0$ and $f\in L^2(\C)$, $\Htp^{s}[f]$ is $C^\infty(\C)$. 
Given a multiindex $\gamma$, there is a constant $C_{|\gamma|}$
so that for $z\in\C$ and 
$R < \min\{R_{\tau p}(z),1\}$ where $R_{\tau p}(z)$ is the constant
from Theorem \ref{thm:sobolev 1},
\[
| X^\gamma \Htp^{s}[f](z) | \leq C_{|\gamma|} R^{-1} s^{- \frac{|\gamma|}2}(1+s^{-1})
 \|f\|_{L^2(\C)}.
\]
\end{lem}

\begin{proof}We can find $\vp_n\in C^\infty_c(\C)$ so that $\vp_n\to f$ in $L^2(\C)$. It follows
immediately from Lemma \ref{lem:X^alpha H^{s} X^beta phi L^2 est} that 
$X^\gamma \Htp^{s}[f] \in L^2(\C)$, and
\[
X^\gamma \Htp^{s}[\vp_n] \to X^\gamma \Htp^{s}[f] 
\]
in $L^2(\C)$, hence
\[
\| X^\gamma \Htp^{s}[f] \|_{L^2(\C)} 
\leq C_{|\gamma|} R^{-1} s^{-\frac{|\gamma|}2} \|f\|_{L^2(\C)}.
\]
From these inequalities, we can show that all 
$\Zstp$ and $\Zbstp$
derivatives of $\Htp^{s}[f]$  are in $L^2(\C)$.
To pass from  $L^2$-bounds of $\Zstp$ and $\Zbstp$ derivatives to a local $L^2$-bound
for ordinary derivatives, we can use the argument of
Theorem \ref{thm:sobolev 1}.  Thus, $\Htp^{s}[f]$ is $C^\infty(\C)$, and 
by Theorem \ref{thm:sobolev 1},
\begin{align*}
\sup_{D(z,R)} |X^\gamma \Htp^{s}[f]|
&\leq \frac{C}{R} \sum_{|\alpha|\leq 2} R^{|\alpha|} \|X^\alpha X^\gamma \Htp^{s}[f] \|_{L^2(\C)}
\leq \frac{C}{R} \sum_{|\alpha|\leq 2} s^{-\frac{|\alpha|+|\gamma|}2} \|f\|_{L^2(\C)} \\
&\leq C_{|\gamma|} R^{-1} s^{-\frac{|\gamma|}2}(1+s^{-1})\|f\|_{L^2(\C)}.
\end{align*}
\end{proof}

\n Recall the following standard fact.
% Lemma: L^2 in each variable for all derivatives means C^\infty
\begin{lem}\label{lem:C^infty together} If $x_1\mapsto \p^\alpha_{x_1} f(x_1,{x_2})$ and 
${x_2}\mapsto \p^\alpha_{x_2} f(x_1,{x_2})$ are in $L^2_{\mathrm{loc}}(\R^n)$ for all multiindices $\alpha$, then
$f\in C^\infty(\R^n\times\R^n)$.
\end{lem}

%\begin{proof} Smoothness is a local property, so we can assume $f$ has compact
%support. Then $\p^{\alpha}_{x_1} f$ and $\p^\beta_{x_2} f$ are in $L^2(\R^n\times\R^n)$ for all
%multiindices $\alpha$ and $\beta$. By the Plancherel
%Theorem, if $\xi$ and $\eta$ are the transform variables
%of $x_1$ and ${x_2}$ respectively, then
%\[
%\|\p^{\alpha}_{x_1}\p^\beta_{x_2} f\|_{L^2(\R^n\times\R^n)} 
%=c \|\xi^\alpha\eta^\beta \hat f\|_{L^2(\R^n\times\R^n)}  
%\leq c \| (|\xi|^{\alpha+\beta} + |\eta|^{\alpha+\beta}) \hat f \|_{L^2(\R^n\times\R^n)}  < \infty.
%\]
%The result follows by the Sobolev Embedding Theorem.
%\end{proof}

We are now ready to prove Theorem \ref{thm:H properties}. We follow the line of argument
in \cite{NaSt00} and show the proof for completeness.
\begin{proof}(Theorem \ref{thm:H properties}) For every multiindex $\alpha$,
Lemma \ref{lem:X H_[s] pointwise est} and Theorem \ref{thm:sobolev 1}
show that the functional on $L^2(\C)$
defined by 
\[
f \mapsto \p^\alpha \Htp^{s}[f]
\]
is bounded. By the Riesz Representation Theorem, a consequence of these facts
is the existence of functions
$\Htp^{\alpha,s,z}(w)$  so that 
\[
\p^\alpha \Htp^{s}[f](z) = \int_{\C} \Htp^{\alpha,s,z}(w) f(w)\, dw.
\]
Define $\Htp^\alpha(s,z,w) = \Htp^{\alpha,s,z}(w)$.
Also set $\Htp(s,z,w) = \Htp^0(s,z,w)$. Then $\Htp^\alpha$ 
is a function on $(0,\infty)\times\C\times\C$ with the property that 
$w\mapsto \Htp^\alpha(s,z,w)$ is in $L^2(\C)$. Thus,
we have
\[
\Htp^{s}[f](z) = \int_{\C} \Htp(s,z,w) f(w)\, dw \label{eqn:H^{s}}
\]
and for every derivative $\p^\alpha_z$,
\[
\p^\alpha_z \left( \int_\C \Htp(s,z,w) f(w)\, dw\right) 
= \int_{\C} \Htp^\alpha(s,z,w) f(w)\, dw.
\]
We will show that $\p^\alpha_z \Htp(s,z,w) = \Htp^{\alpha}(s,z,w)$. 
Let $\vp,\psi \in C^\infty(\C)$. By the Schwartz Kernel Theorem,
\begin{align*}
\langle \p^\alpha_z \Htp^{s}[\psi],\vp\rangle_{\C}
&= \langle (-1)^{|\alpha|} \Htp^{s}[\psi], \p^\alpha_z\vp \rangle_{\C} 
= \langle (-1)^{|\alpha|} \Htp^{s}, \psi\otimes \p^\alpha_z\vp\rangle_{\C\times\C} \\
&= \langle \p_z^\alpha \Htp^{s}, \psi\otimes\vp\rangle_{\C\times\C}
= \langle (\p^\alpha_z \Htp^{s})[\psi],\vp\rangle_{\C}.
\end{align*}
Thus, we have shown that
\begin{equation}
\p^\alpha_z \Htp(s,z,w) = \Htp^\alpha(s,z,w) \label{eqn:H alpha deriv}
\end{equation}
in $\mathcal{D}'(\C)$ and that
for each $\alpha$, $w\mapsto \p^\alpha_z \Htp(s,z,w)$ is an $L^2(\C)$ function.

Next, we know that $\Htp^{s}$ is self-adjoint, so
\[
\int_{\C} \Htp^{s}[\psi](z)\overline{\vp(z)} \, dz = \int_{\C} \psi(w)
\overline{\Htp^{s}[\vp](w)}\, dw.
\]
As an immediate consequence of this equality and \eqref{eqn:H^{s}},
we have
\[
\int_\C\int_\C \Htp(s,z,w)\psi(w)\overline{\vp(z)} \, dw dz 
= \int_\C\int_\C \overline{\Htp(s,w,z)} \psi(w)\overline{\vp(z)}\, dz dw. 
\]
It follows that $\Htp(s,z,w) = \overline{\Htp(s,w,z)}$, conclusion  (a).

As a consequence of  (a)  and the fact that $w \mapsto \Htp^{s}(s,z,w)$
belongs to $L^2(\C)$, $z\mapsto \Htp^{s}(s,z,w)$ 
belongs to $L^2(\C)$. By equations \eqref{eqn:H alpha deriv}, 
it follows that every $z$ derivative also belongs to $L^2$. Thus by
Lemma \ref{lem:C^infty together}, $\Htp(s,z,w)$ is $C^\infty(\C\times\C)$ for 
fixed $s>0$.

We know $\Boxtp^j \Htp^{s} = \Htp^{s} \Boxtp^j$. 
The implication of the self-adjointness of 
$\Boxtp$ is that on the kernel side, 
\[
\Boxtpz^j \Htp(s,z,w) = (\Boxtpw\sh)^j \Htp(s,z,w).
\] 
From this,  (c)  follows
quickly because 
$\Boxtpz^{j+k} \Htp^{s} = \Boxtpz^j \Htp^{s} \Boxtpw^k$. 

Next, by Theorem \ref{thm:semigroups} (e), 
$\left(\frac{\p}{\p s} + \Boxtp\right)\big[\Htp^{s}[f]\big](z) =0$. 
Fixing $z\in\C$, integration 
against test functions in $(0,\infty)\times\C$ shows that in 
$\mathcal{D}'\big((0,\infty)\times\C\big)$,
\begin{align*}
0 &= \left\langle \left(\frac{\p}{\p s} + \Boxtp\right)
\big[\Htp^{s}[f]\big](z), \vp\right\rangle \\
&= \iint_{(0,\infty)\times\C} \Htp(s,z,w)\left(-\frac{\p}{\p s} + \Boxtpz\sh\right)\vp(s,z)f(w) \, dw  ds \\
&= \iint_{(0,\infty)\times\C}\left(\frac{\p}{\p s} + \Boxtpz\right) 
\Htp(s,z,w)\vp(s,z) f(w) \, dw  ds.
\end{align*}
Thus, 
$\frac{\p}{\p s} \Htp(s,z,w) = -\Boxtpz \Htp(s,z,w)$ .
Then we have 
\begin{align*}
\frac{\p^2}{\p s^2} \Htp(s,z,w) &= -\frac{\p}{\p s}\Boxtpz \Htp(s,z,w)\\ 
&= -\Boxtpz \frac{\p}{\p s} \Htp(s,z,w) = \Boxtpz^2 \Htp(s,z,w).
\end{align*}
Iterating this argument shows
\begin{equation}
\frac{\p^j}{\p s^j} \Htp(s,z,w)= (-1)^j\Boxtpz^j \Htp(s,z,w)  \label{eqn:H diffeq}.
\end{equation}
We know, however, that for each $j$, $\Boxtpz^j\Htp(s,z,w)\in 
L^2_{\mathrm{loc}}\big((0,\infty)\times\C\times\C\big)$. As before,
this is enough to show
$\Htp \in C^\infty\big((0,\infty)\times\C\times\C\big)$. In particular,  \eqref{eqn:H diffeq}
hold in the classical sense, so (b) is proved.

For $\alpha$, $\beta$, and $j$, Lemma \ref{lem:X^alpha H^{s} X^beta phi L^2 est}
shows that there is a constant $C_{\alpha,\beta,j}$ so that for $\vp\in C^\infty_{c}(\C)$,
\[
\Big \| X^\alpha \Boxtp^j \Htp^{s} \big[X^\beta [\vp]\big] \Big\|_{L^2(\C)}
\leq C_{\alpha,\beta} s^{-\frac{|\alpha|+|\beta|}2-j} \|\vp\|_{L^2(\C)}.
\]
Then by Theorem \ref{thm:sobolev 1}, for $R<R_{\tau p}(z)$,
\begin{align*}
\sup_{D(z,R)} \big| X^\alpha \Boxtp^j \Htp^{s}\big[X^\beta [\vp]\big] \big|
&\leq \frac CR \sum_{|\gamma|\leq 2} R^{|\gamma|}
\big \| X^{\alpha+\gamma} \Boxtp^j \Htp^{s} \big[X^\beta [\vp]\big]\big \|_{L^2(\C)} \\
&\leq \frac CR s^{-\frac{\alpha+\beta}2-j}(1+ s^{-1}) \|\vp\|_{L^2(\C)}.
\end{align*}

\n Also, since $\Htp(s,z,w)\in C^\infty\big((0,\infty)\times\C\times\C\big)$, 
\begin{align*}
X^\alpha \Boxtp^j \Htp^{s}[X^\beta\vp](z) 
&= \int_{\C} X_z^\alpha \Boxtp^j \Htp(s,z,w) X_w^\beta\vp(w) \, dw \\
%&= \int_{\C} X_z^\alpha U_w^\beta \Boxtp^j \Htp(s,z,w) \vp(w) \, dw\\
&= (-1)^{j+|\beta|} \int_{\C} \frac{\p^j}{\p s^j} X_z^\alpha U_w^\beta \Htp(s,z,w)\vp(w)\, dw.
\end{align*}
From the reverse H\" older inequality and our previous estimate,
\[
\left(\int_\C \left| \frac{\p^j}{\p s^j} X_z^\alpha U_w^\beta \Htp(s,z,w)\right|^2\, dw \right)^{\frac 12}
\leq \frac {C_{\alpha,\beta, j}}R s^{-\frac{\alpha+\beta}2-j}(1+ s^{-1}).
\]
This is 
(d) of the Theorem. From (a),
we can interchange the roles of $z$ and $w$ to prove (e). This proves
the theorem.
\end{proof}

%%%%%%%%%%%%%%%%%%%%
%
% SECTION: FUNDAMENTAL SOLUTION FOR H
%
%%%%%%%%%%%%%%%%%%%%
\section{A Fundamental Solutions for $\H$ on $\R\times\C$} \label{sec:fundemantal solutions}

%subsection: Fundamental Solutions for H
\subsection{A Fundamental Solution for $\H$ and a Relative Fundamental Solution for $\Hw$
on $\R\times\C$} \label{subsec:fundamental solutions}

Define distributions $\Htp^{z}$, $\Hwtp^{z}$  on $\R\times\C$ by
\begin{defn}\label{defn:H_z}
For $\psi\in C^\infty_c(\R\times\C)$, set
\[
\langle \Htp^{z},\psi \rangle = \lim_{\ep\to 0} \int_{\ep}^\infty \int_{\C} \Htp(s,z,w) \psi(s,w)\, dw ds.
\]
\end{defn}

The kernel of this distribution is
\[
\Htp^{z}(s,w) = \begin{cases} \Htp(s,z,w) & \text{if } s >0 \\
 0 &\text{if } s\leq 0. \end{cases}
\]

We need to prove the limit defining $\Htp^z$ exists and defines a distribution
on $\R\times\C$. To do so, we need the following lemma which 
gives control over pointwise bounds on $|e^{-s\Boxtp}f-f|$.
% Lemma: POINTWISE BOUNDS ON |F - E^-S\BOX F
\begin{lem}\label{lem:pointwise bounds}
There is a constant depending only on the degree of $p$ so that if $f\in \Dom(\Boxtp^j)$
for $j \leq 2$, then for any $z\in\C$ and $0<R<R_{\tau p}(z)$ and $z\in\C$,
\[
\sup_{D(z,R)} \big| f(w) - e^{-s\Boxtp}[f](w)\big| 
\leq C \frac sR \left(\|\Boxtp  [f] \|_{L^2(\C)} + R^{2} \|\Boxtp^{2}[f] \|_{L^2(\C)}\right).
\]
\end{lem}

\begin{proof} Let $f\in \Dom(\Boxtp^j)$, $0 \leq j \leq 2$. 
By Theorem \ref{thm:sobolev 1} and Theorem \ref{thm:semigroups} (c), we have
\begin{align*}
&\sup_{D(z,R)} \big| f(w) - e^{-s\Boxtp}[f](w)\big| 
\leq C \frac 1R \sum_{j=0}^{1} R^{2j} \big \| \Boxtp^j[f] 
- \Boxtp^j[e^{-s\Boxtp}f] \big\|_{L^2(\C)} \\
&=\frac CR \sum_{j=0}^{1} R^{2j} \big\|(I-e^{-s\Boxtp})[\Boxtp^j f]\big\|_{L^2(\C)} 
\leq C \frac sR\sum_{j=0}^{1} R^{2j} \|\Boxtp^{j+1}[f] \|_{L^2(\C)}.
\end{align*}
\end{proof}

\begin{lem}\label{lem:H_z defines  a distribution}
For each $z\in \C$, the limit defining  $\Htp^{z}$  exists and defines a distribution on 
$\R\times\C$. 
\end{lem}

\begin{proof} Let $\psi\in C^\infty_c(\R\times\C)$. Then there is a closed, bounded interval 
$I\subset \R$ and a compact set $K\subset\C$ so that $\supp\psi\in I\times K$. Set
$\psi_s(z) = \psi(s,z)$. Then $ \{\psi_s\}\subset C^\infty_c(\R)$ with each element having
support in $K$. If $0<\epsilon_1<\epsilon_2$, then
\begin{align*}
\int_{\ep_1}^{\ep_2} \Htp(s,z,w)&\psi(s,w)\, dw\, ds 
= \int_{\ep_1}^{\ep_2} e^{-s\Boxtp}[\psi_s](z)\, ds \\
&= \int_{\ep_1}^{\ep_2} e^{-s\Boxtp}[\psi_s](z) - \psi(s,z)\, ds + \int_{\ep_1}^{\ep_2}\psi(s,z)\, ds.
\end{align*}
From Lemma \ref{lem:pointwise bounds} and H\" older's inequality, we have (with $R<R_{\tau p}(z)$),
\begin{align*}
\Big| \int_{\ep_1}^{\ep_2} \int_{\C} &\Htp(s,z,w) \psi(s,w)\, dw ds \Big| \\
&\leq C \frac {\ep_2}R \sum_{j=0}^{1} R^{2j}
\int_{0}^{\ep_2} \|\Boxtp^{j+1}[\psi_s]\|_{L^2(\C)}\, ds + \ep_2 \|\psi\|_{L^\infty(\R\times\C)}\\
&\leq C\frac{\ep_2^{\frac 32}}R \sum_{j=0}^{1}\|\Boxtp^{j+1}[\psi]\|_{L^2(\R\times\C)}
+ \ep_2\|\psi\|_{L^\infty(\R\times\C)}.
\end{align*}
These last terms go to 0 as $\ep_2\to 0$, so the limit defining $\Htp^{z}$ exist. \end{proof}

% THEOREM: H_Z IS THE FUNDAMENTAL SOLUTION
\begin{thm}\label{thm:fundamental solution}
In $\mathcal{D'}(\R\times\C)$, 
\[
(\p_s + \Boxtpw\sh)(\Htp^{z}) = \delta_0\otimes \delta_z.
\]
\end{thm}

\begin{proof} Let $\psi\in C^\infty_c(\R\times\C)$. Then
\begin{align*}
\langle(\p_s + \Boxtpw\sh)(\Htp^{z}),\psi\rangle 
=& \langle \Htp^{z}, (-\p_s + \Boxtpw)\psi\rangle\\
= &-\lim_{\ep\to0} \int_{\ep}^{\infty}\int_\C \Htp(s,z,w)\p_s\psi(s,w)\, dw ds\\
&+ \lim_{\ep\to 0}\int_{\ep}^\infty \int_\C \Htp(s,z,w) \Boxtpw\psi(w,s)\, dw ds.
\end{align*}
Since $s$ is bounded away from $0$ and $\Htp\in C^\infty\big((0,\infty)\times\C\times\C\big)$, 
the first term yields
\begin{align*}
&-\int_{\ep}^{\infty}\int_\C \Htp(s,z,w)\frac{\p\psi}{\p s}(s,w)\, dw ds\\
&= - \int_{\ep}^{\infty}\frac{\p}{\p s}\int_\C \Htp(s,z,w)\psi(s,w)\, dw ds 
+\int_{\ep}^{\infty}\int_\C \frac{\p}{\p s}\Htp(s,z,w)\psi(s,w)\, dw ds\\
&= \int_{\C} \Htp(\ep,z,w)\psi(\ep,w)\, dw + \int_{\ep}^\infty \int_{\C} 
\frac{\p}{\p s} \Htp(s,z,w)\psi(s,w)\, dw ds.
\end{align*}
Also,
\[
\int_{\ep}^\infty\int_{\C} \Htp(s,z,w) \Boxtpw\psi(s,w)\, dw ds
= \int_{\ep}^\infty \int_{\C} \Boxtpw\sh \Htp(s,z,w)\psi(s,w)\, dw ds.
\]
Using Theorem \ref{thm:H properties} (b) and adding our equalities together, we have
\begin{align*}
&\int_{\ep}^\infty \int_\C-\Htp(s,z,w)\p_s\psi(s,w) + \Htp(s,z,w) \Boxtpw\psi(s,w)\, dw ds\\
&= \int_{\C}\Htp(\ep,z,w)\psi(\ep,w)\, dw + \int_{\ep}^\infty\int^\infty_{\C} (\p_s+\Boxtpw\sh)\Htp(s,z,w)
\psi(s,w)\, dw ds \\
&= \int_{\C} \Htp(\ep,z,w)\psi(\ep,w)\, dw.
\end{align*}
Hence
\begin{align*}
 \langle (\p_s+\Boxtpw\sh)[\Htp^{z}],\psi\rangle &= \lim_{\ep\to 0} \int_{\C} \Htp(\ep,z,w)\psi(\ep,w)\, dw\\
&= \lim_{\ep\to 0} e^{-\ep\Boxtp}[\psi_\ep](z) = \psi(0,z) = \langle \delta_0\otimes\delta_z,\psi\rangle.
\end{align*}
\end{proof}

\section{Estimates on $\frac{\p^n}{\p s^n}Y^\alpha\Htp(s,z,w)$}\label{sec:heat estimate}
In this section, we prove Theorem \ref{thm:bounds on H}, the
result on pointwise estimates of $|X_z^I U_w^J \Htp(s,z,w)|$.
We begin the section with a study of how the heat kernel behaves under scaling.

%%%%%%%%%%%%%%%
%
%		SECTION: SCALING ESTIMATES
%
%%%%%%%%%%%%%%%
\subsection{Scaling and the Heat Kernel}\label{subsec:scaling}
The structure of $\Boxtp$ is critical in this section. Expanding $\Boxtp$, we have
\begin{align}
&\Boxtp = -\left(\frac{\p}{\p\z} + \tau\frac{\p p}{\p \z}\right)
\left( \frac{\p }{\p z} - \tau\frac{\p p}{\p z} \right) \nonumber\\
&= -\frac{\p^2 }{\p z\p\z} + \tau \frac{\p^2 p}{\p z\p\z}
+ \tau^2 \frac{\p p}{\p z}\frac{\p p}{\p\z}
+\tau\left( \frac{\p p}{\p z} \frac{\p }{\p\z} - \frac{\p p}{\p\z}\frac{\p }{\p z}\right) \label{eqn:Boxz}\\
&= -\frac14 \triangle + \frac14 \tau \triangle p+ \frac {\tau^2}4|\nabla p|^2  + \frac i2\tau 
\left(\frac{\p p}{\p x_1}\frac{\p}{\p x_2} - \frac{\p p}{\p x_2}\frac{\p}{\p x_1}\right) \label{eqn:Boxx}
\end{align}

Let $p_0(w) = p(w)$ and fix $z_0\in\C$. Let $p_1(w) = p_0(w+z_0)$. Our first scaling result is:
% Proposition: translation of the heat kernel
\begin{prop}\label{prop:translation} 
\[
H_{\tau p_0}(s,z+z_0,w+z_0) = H_{\tau p_1}(s,z,w).
\]
\end{prop}

\begin{proof} Fix $z_0\in\C$.
Let $A_{z_0}[f](z) = f(z-z_0)$. $A_{z_0}$ is an isometry on $L^2(\C)$, and 
\begin{align*}
\HH_{\tau p_1}&[f](z)
= -\frac 14\triangle f(z) + \frac{\tau}4\triangle p_0(z+z_0)f(z) 
+ \frac{\tau^2}4 |\nabla p_0(z+z_0)|^2 f(z) \\
&\ \ + \frac i2 \tau\left(\frac{\p p_0}{\p x_1}(z+z_0) 
\frac{\p f}{\p x_2}(z) - \frac{\p p_0}{\p x_2}(z+z_0)  \frac{\p f}{\p x_1}(z)\right) \\
&= A_{z_0}^{-1} \Bigg[-\frac 14\triangle f(z-z_0) + \frac{\tau}4\triangle p_0(z)f(z-z_0) 
+ \frac{\tau^2}4 |\nabla p_0(z)|^2 f(z-z_0) \\
&\ \ + \frac i2 \tau\left(\frac{\p p_0}{\p x_1}(z) 
\frac{\p f}{\p x_2}(z-z_0) - \frac{\p p_0}{\p x_2}(z)  \frac{\p f}{\p x_1}(z-z_0)\right) \Bigg]\\
&= A_{z_0}^{-1} \HH_{\tau p_0} A_{z_0}[f](z).
\end{align*}
Also, if $\psi\in C^\infty_c(\C\times\R)$,
\[
A_{z_0}^{-1} (\delta_0 \otimes\delta_z) A_{z_0}\psi(s,w)
= A_{z_0}^{-1}\psi(0,z-z_0) = \psi(0,z),
\]
and
\begin{align*}
A_{z_0}^{-1} (\delta_0 &\otimes\delta_z) A_{z_0}\psi(s,w)
= A_{z_0}^{-1} \HH_{\tau p_0} \int_{0}^\infty \int_{\C} H_{\tau p_0}(s,z,w) \psi(s,w-z_0)\, dw ds\\
&= A_{z_0}^{-1}\int_{0}^\infty \int_{\C} \Big(\frac{\p}{\p s} + \Box_{\tau p_0,z}\Big) 
H_{\tau p_0}(s,z,w) \psi(s,w-z_0)\, dw ds\\
&= A_{z_0}^{-1}\int_{0}^\infty \int_{\C} \Big(\frac{\p}{\p s} + \Box_{\tau p_0,w}\sh\Big)
H_{\tau p_0}(s,z,w) \psi(s,w-z_0)\, dw ds\\
&= A_{z_0}^{-1}\int_{0}^\infty \int_{\C} \Big(\frac{\p}{\p s} +  \Box_{\tau p_1,w}\sh\Big)
H_{\tau p_0}(s,z,w+z_0) \psi(s,w-z_0)\, dw ds\\
&= \int_{0}^\infty \int_{\C} \Big(\frac{\p}{\p s} +  \Box_{\tau p_1,w}\sh\Big)
H_{\tau p_0}(s,z+z_0,w+z_0) \psi(s,w-z_0)\, dw ds.
\end{align*}
Thus, $H_{\tau p_0}(s,z+z_0,w+z_0) = H_{\tau p_1}(s,z,w)$.
\end{proof}

The method of proof in Proposition \ref{prop:translation} also proves the following two scaling results.
%Proposition: rotation of the heat kernel
\begin{prop}\label{prop:rotation}If $z_0\in\C$ and
\[
p_2^{z_0}(w) = \sum_{j,k\geq 1} \A{jk}{z}(w-z_0)^j\overline{(w-z_0)}^k,
\]
then
\[
H_{\tau p_2^{z_0}}(s,z,w) = e^{i \tau  T(w,z_0)}H_{\tau p_1}(s,z_0,w).
\]
\end{prop}

Let $T_\lam \psi(s,w) = \lam^2 \psi(\lam^2 s, \lam w)$ and
$T_\lam f(w) = \lam f(\lam w)$ act on functions on $\R\times\C$ and $\C$, respectively. In either
case, $T_\lam$ is an isometry on $L^2$. Our final proposition in this section 
investigates
conjugating $\H$ by $T_\lam$. Let $\psi_\lam(s,w) = \psi(\lam^2 s,\lam w)$
and $f_\lam(w) = f(\lam w)$. 
%Proposition: dilation of p
\begin{prop} \label{prop:dilation} If $p_3^\lam = p_2(z/\lam)$, then
\[
\frac 1{\lam^2} H_{\tau p_2^{z_0}}(s/\lam^2, z/\lam, w/\lam) 
= H_{\tau p_3^{\lam}}(s,z,w).
\]
\end{prop}

%%%%%%%%%%%%%%%%
%
%	SECTION: ESTIMATE ON H(S,Z,W)
%
%%%%%%%%%%%%%%%%%
\subsection{Pointwise Estimates of $|\Htp(s,z,w)|$}\label{subsec:est on H}

We first show 
that $|\Htp(s,z,w)|$ has Gaussian decay. To do so, we will find it convenient to work
in real variable notation instead of complex notation.
As such, let $x = (x_1,x_2)$ and $y = (y_1, y_2)$. Our first goal is to
prove:
% Theorem: Gaussian Decay for H(s,x,y)
\begin{thm}\label{thm:gaussian} If $e^{-s\Boxtp}[f](x) = \int_{\C} \Htp(s,x,y)f(y)\, dy$,
then the heat kernel $\Htp(s,x,y)$, satisfies the estimate
\[
|\Htp(s,x,y)| \leq \frac{1}{\pi s} e^{-\frac{|x-y|^2}s}.
\]
\end{thm}

\begin{proof} We will use the Feynman-Kac-It\^{o} formula from \cite{Simon79}.
Let $dx$ be Lebesgue
measure on $\R^2$ and let $(B,\mathfrak{B},dP)$ be a measure space of sample paths
for a 2-dimensional Brownian motion $b(s)$. Let $d\mu = dP\otimes dx$ be Wiener measure on
$B\times\R^2$ and let $\omega(s) = x+b(s)$. If we let
\[
a(x) = \tau \left(  -\frac{\p p}{\p x_2}, \frac{\p p}{\p x_1} \right)
\qquad\text{and}\qquad
V(x) = \frac \tau2 \triangle p(x),
\]
then for $f\in C^2(\R^2)$,
\[
\frac12(-i\nabla - a)^2f + Vf
= -\frac 12\triangle f+ \frac i2(\nabla\cdot a)f + ia\cdot\nabla f + \frac 12 |a|^2 f + Vf.
\]
But  $\nabla\cdot a = \tau(-\frac{\p^2 p}{\p x_1\p x_2} + \frac{\p^2 p}{\p x_1\p x_2})=0$
and
$\frac 12 |a|^2 = \frac 12\tau^2|\nabla p|^2$. Thus,
\[
\frac12(-i\nabla - a)^2f + Vf = 2\Box_{\tau p},
\]
so $2\Boxtp$ is the quantum mechanical energy operator for a particle in a magnetic
field with vector potential $a(x)$ and electric potential $V$. The Feynman-Kac-It\^{o} formula for
$f,g\in C^\infty_c(\R^2)$ is
\begin{equation}\label{eqn:FKI}
\Big( e^{-2s\Boxtp} f,g\Big) = \int e^{F(s,\omega)}f(\omega(s))\overline{g(\omega(0))}\,d\mu
\end{equation}
where 
\[
F(s,\omega) = -i\int_0^s a(\omega(t))\cdot d\omega(t) - \frac i2 \int_{0}^s (\nabla \cdot a)(\omega(t))\,dt
- \int_0^s V(\omega(t))\, dt.
\]
$b(s)$ has 2-dimensional normal distribution
with covariance $s$, so we can rewrite \eqref{eqn:FKI} as follows:
\begin{align*}
\iint_{\R^2\times\R^2} &\Htp(2s,x,y) f(y)\overline{g(x)}\, dydx
=  \int e^{F(s,\omega)}f(\omega(s))\overline{g(\omega(0))}\,d\mu\\
&= \int_{\R^2} E\big[ e^{F(s,\omega)}f(\omega(s)) \overline{g(\omega(0))}\big]\,  dx \\
&= \int_{\R^2} E\big[E[e^{F(s,\omega)}f(\omega(s)) \overline{g(\omega(0))}\big|
\omega(0)=x,\omega(s)=x+y]\big]\, dx\\
&=\int_{\R^2} E\big[E\big[e^{F(s,\omega)}f(\omega(s))\big|
\omega(0)=x,\omega(s)=x+y\big]\big] \overline{g(x)}\, dx\\
&= \frac 1{2\pi s}\iint_{\R^2\times\R^2} e^{\tilde F(s,x+y)}f(x+y)\overline{g(x)} e^{-\frac {|y|^2}{2s}}\, dydx\\
&= \frac 1{2\pi s}\iint_{\R^2\times\R^2}e^{\tilde F(s,y)}f(y)\overline{g(x)} e^{-\frac {|x-y|^2}{2s}}\, dydx.
\end{align*}
Thus, $\Htp(2s,x,y) = \frac 1{2\pi s}e^{\tilde F(s,y)}e^{-\frac {|x-y|^2}{2s}}$
for some $\tilde F(s,y)$ satisfying $|e^{\tilde F(s,y)}|\leq 1$. 
%That $|e^{\tilde F(s,y)}|\leq 1$ follows from $|e^{ F(s,\omega)}|\leq 1$.
\end{proof}
A critically important fact about the Feynman-Kac-It\^{o} formula is the requirement that
$V \geq 0$. When $\tau<0$, $V\leq 0$, and the argument from 
Theorem \ref{thm:gaussian} fails. Even if we could use the argument, the real part of
$e^{F(s,\omega)}$ is $e^{-\int_0^s V(\omega(t))\, dt}$, a term that we would expect to be
very large. In fact, when $\tau<0$ $\Htp(s,z,w)$ only satisfies Gaussian decay near
$s=0$ \cite{Rai06}, so an analog to Theorem \ref{thm:gaussian} is false.

We now turn to proving a large time decay estimate for $\Htp(s,z,w)$.
Let $\mathfrak{P}$ be the set of polynomials  of $\deg(p)$ whose 
coefficients (in absolute value) sum to 1. 
We can identify the set of polynomials of $\deg(p)$ with
$\R^n$ for some $n$, and under this identification,
$\mathfrak{P}$ is identified with the unit sphere, a compact set. 
Having constants depending only on $\mathfrak{P}$ is essential for estimates
obtained through scaling.

%
%	THEOREM: TIME DECAY ESTIMATE ON H(s,z,w)
%
\begin{thm}\label{thm:heat est}  If $e^{-s\Boxtp}[f](z) = \int_{\C} \Htp(s,z,w)f(w)\, dw$,
then there exist constants $C_1$ and $C_2$ which
depend on the degree of $p$ so that
\[
|\Htp(s,z,w)| 
\leq \frac{C_1}s e^{-C_2\frac{s}{\mu(z,1/\tau)^2}}e^{-C_2\frac{s}{\mu(w,1/\tau)^2}}.
\]
\end{thm}

\begin{proof} By the conjugate symmetry of $\Htp$, it is enough to show the bound
$\frac{C_1}s e^{-C_2\frac{s}{\mu(z,1/\tau)^2}}$.
From Proposition \ref{prop:translation}, there exists
a polynomial $p_1(w) = p_0(w+z)$ so that
\[
\Htp(s,z,w) = H_{\tau p_1}(s,0,w-z),
\]
so we can reduce to the case of estimating $H_{\tau p_1}(s,0,w)$. By
Proposition \ref{prop:rotation}, for
\[
p_2(w) = \sjk \A{jk}{0}w^j\bar w^k,
\]
we have
\[
e^{i \tau T(w,0)} H_{\tau p_1}(s,0,w) = H_{\tau p_2}(s,0,w),
\]
so it is enough to estimate $|H_{\tau p_2}(s,0,w)|$. Observe that $p_2(w)$ has the
property $\frac{\p^k p_2}{\p z^k}(0) = \frac{\p^k p_2}{\p\z^k}(0)=0$ 
for all $k$. If we set $\lam = \mu(z,1/\tau)^{-1}$ and
$p_3(w) = p_2(\frac w\lam)$, then $p_3 \in \mathfrak{P}$ since
\[
\tau \sjk \frac{1}{j!k!} \left|\frac{\p^{j+k}p_2}{\p z^j \p\z^k}(0)\right|
z^j\z^k \lam^{-j-k} = \tau\Lambda\big(0,\mu(0,\tfrac 1\tau)\big) \sim 1.
\]
From Proposition
\ref{prop:dilation}, 
\[
\frac{1}{\lam^2} H_{\tau p_2}(\tfrac s{\lam^2},0, \tfrac w\lam) 
= H_{p_3}(s,0,w).
\]
We now estimate $|H_{p_3}(s,w,0)|$. Let $h(s,w) = H_{p_3}(s,w,0)$.  
By Theorem \ref{thm:H properties} (b), 
$\frac{\p h}{\p s} = \overline{Z}_{p_3} Z_{p_3} h$.
Let $g(s) = \int_{\C} |h(s,w)|^2\, dw$.
From \cite{Christ91}, there exists $C = C_\mathfrak{P}$ so that
$\|f\|_{L^2(\C)}\leq C\|Z_{p_3} f\|_{L^2(\C)}$. Thus
\begin{align*}
g'(s) &= \int_{\C} \frac{d}{ds}\big(h(s,w)\overline{h(s,w)}\big)\, dw
= 2\Rre \int_{\C} \frac{\p h}{\p s}(s,w)\overline{h(s,w)}\, dw \\
&= 2\Rre \int_{\C} \overline{Z}_{p_3} Z_{p_3} 
h(s,w)\overline{h(s,w)}\, dw 
= -2 \int_{\C}\big| Z_{p_3} h(s,w)\big|^2\, dw \\
&\leq - C\int_{\C}\big|  h(s,w)\big|^2\, dw = -Cg(s).
\end{align*}
Since $g(s)>0$, $\frac{g'(s)}{g(s)} \leq -C$, and integrating from $\frac s2$ to $s$,
we have
\[
g(s) \leq g(\tfrac s2) e^{-C s} \leq C_1 \frac {e^{-C s}}s,
\]
where the last inequality follows from Theorem \ref{thm:gaussian}. The
constant $C_1$ does not depend on $p_3$ (or $\mathfrak{P}$).

Next, $e^{-s\Box_{p_3}}$ is a semigroup, so 
$e^{-s\Box_{p_3}}e^{-s\Box_{p_3}}f(z) 
= e^{-2s\Box_{p_3}}f(z)$. On the kernel side, this means 
we have the reproducing identity
\[
H_{p_3}(2s,z,w) = \int_{\C} H_{p_3}(s,z,v)
H_{p_3}(s,v,w)\, dv,
\]
and an application of Cauchy-Schwarz yields
\begin{align*}
\big| H_{p_3}(2s,0,w)\big|
&\leq \left(\int_{\C} \big|H_{p_3}(s,0,v)\big|^2\, 
dv\right)^{\frac 12} \left(\int_{\C} \big|H_{p_3}(s,v,w)\big|^2\, 
dv\right)^{\frac 12} \leq C_1\frac {e^{-C s}}s.
\end{align*}
Undoing the scaling finishes the proof.
\end{proof}

The motivation for using $g'(s)$ and the reproducing
identity was \cite{Fabes93}.

%%%%%%%%%%%%%%%%%%
%
%	SECTION: DERIVATIVE ESTIMATES
%
%%%%%%%%%%%%%%%%%%%
\subsection{Derivative Estimates}\label{subsec:derivative}
The derivative estimates are proven in a series of lemmas. The most accessible case
is proven first and each successive lemma builds on the previous calculation. Each
$L^2$ estimate at one step is used to prove a pointwise estimate in the next.
Define the decay term $D(s,x,y)$ to be
\begin{equation}\label{eq:D(s,x,y)}
D(s,x,y) = e^{-\frac{|x-y|^2}{2 s}} e^{-C_2\frac{s}{\mu(x,1/\tau)^2}}
e^{-C_2\frac{s}{\mu(y,1/\tau)^2}}
\end{equation}
where $C_2$ is the constant from Theorem \ref{thm:heat est}. Also, let
\[
I_r(s) = (s-r^2,s) \text{ and }Q_r(s,x) = I_r(s)\times D(x,r).
\]

We need  a version of the  subsolution estimate from 
\cite{Kurata00}. 
%
% LEMMA: SUBSOLUTION ESTIMATE
\begin{lem}\label{lem:kur sub} If $(s_0,z_0)\in (0,\infty)\times \C$ and
$u(s,z)$ is a $C^2$ solution of
\[
\frac{\p u}{\p s} + \Boxtp u =0
\]
on $Q_{2r}(s_0,z_0)$. Then if $\tau>0$, there exists $C>0$ so that
\[
\sup_{(s,z)\in Q_{r/2}(s_0,z_0)} |u(s,z)|
\leq \frac C{r^2} \iint_{Q_{2r/3}(s_0,z_0)} |u(s,z)|^2\, dzds.
\]
\end{lem}

% Proposition: L^2 estimates on \frac{\p}{\p s^n} H(s,x,y)
\begin{prop}\label{prop:L^2 s deriv} There exists $C_n$ so that for $0 < r < \frac {\sqrt {s_0}}{16}$,
\[
\left \| \frac{\p^n \Htp}{\p s^n}(\cdot,x,\cdot) \right\|_{L^2(Q_r(s_0,y_0))} \leq \frac C{ s_0^n}.
\]
\end{prop}

\begin{proof}We have

\begin{align}
\left\| \frac{\p^n \Htp}{\p s^n}(\cdot, x_0, \cdot)\right\|_{L^2(Q_r(s_0,y_0))}^2
&= \int_{I_r(s_0)} \left|\left( \int_{D(y_0,r)} \left|\frac{\p^n \Htp}{\p s^n}(s, x_0,y)\right|^2\, dy \right)^{1/2}
\right|^2\, ds\nn\\
&= \int_{I_r(s_0)}\bigg| \sup_{\atopp{\vp\in C^\infty_{c}(D(y_0,r))}{\|\vp\|_{L^2}=1}} \int
\frac{\p^n \Htp}{\p s^n}(s, x_0,y)\vp(y)\, dy\bigg|^2\, ds \nn\\
&=  \int_{I_r(s_0)}\bigg| \sup_{\atopp{\vp\in C^\infty_{c}(D(y_0,r))}{\|\vp\|_{L^2}=1}} 
\frac{\p^n }{\p s^n}\Htp^{s}[\vp](x_0)\bigg|^2\, ds. \label{eq:crazy L2 sup}
\end{align}
The key to the proof is that $\frac{\p^n }{\p s^n}\Htp^{s}[\vp](x)$ satisfies 
$( \frac{\p}{\p s} + \Boxtpx)\frac{\p^n }{\p s^n}\Htp^{s}[\vp](x)=0$. 
By Lemma \ref{lem:kur sub} 
and Theorem \ref{thm:semigroups} (d), estimating an arbitrary term from the supremum in
\eqref{eq:crazy L2 sup} yields
\begin{align}
\left| \frac{\p^n }{\p s^n}\Htp^{s}[\vp](x_0)\right| 
&\leq \frac{C}{r^2} \left( \iint_{Q_{r}(s,x_0)}\left| \frac{\p^n}{\p t^n}
\Htp^{t}[\vp](x)\right|^2\, dx dt\right)^{1/2}\nn\\
&\leq \frac{C}{r^2} \left( \int_{I_{\sqrt 2r}(s_0)}
\left\| \frac{\p^n}{\p t^n}\Htp^{t}[\vp]\right\|_{L^2(D(x_0,r))}^2 dt\right)^{1/2}\nn\\
&\leq \frac{C}{r^2}\left( \int_{s_0-2r^2}^{s_0} \frac{1}{t^{2n}}\, dt\right)^{1/2} 
\leq \frac{C}{r s_0^{n}}. \label{eq:crazy L2 sup est}
\end{align}
Putting \eqref{eq:crazy L2 sup} into \eqref{eq:crazy L2 sup est}, we have
\[
\left\| \frac{\p^n \Htp}{\p s^n}(\cdot, x_0, \cdot)\right\|_{L^2(Q_r(s_0,y_0))}
\leq C\left(\int_{I_r(s_0)} \frac{1}{r^2 s_0^{2n}}\, ds\right)^{1/2} =  \frac{C}{s_0^n}.
\]
\end{proof}

%% Lemma: \frac{\p}{\p s^n} H(s,x,y)
\begin{lem}\label{lem:s derivatives}
Let $n_1$, $n_2$, $n_3\geq 0$ and $n=n_1+n_2+n_3$. Then there exists $C_n>0$ so that
\[
\left|\frac{\p^{n_1}}{\p s^{n_1}} \Boxtpx^{n_2} (\Boxtpy\sh) ^{n_3}\Htp(s,x,y)\right | 
\leq \frac{C_n}{s_0^{n + 1}} D(s,x,y)^\frac{1}{2}.
\]
\end{lem}

\begin{proof} Since $\Htp$ satisfies $( \frac{\p}{\p s} + \Boxtpx)\Htp(s,x,y)=0$ when $s\neq 0$ or $x\neq y$,
it is enough to show the estimate for $H_n(s,x,y) = \frac{\p^n}{\p s^n}\Htp(s,x,y)$. Proof by induction. 
The base case follows from combining Theorem \ref{thm:gaussian} and 
Theorem \ref{thm:heat est}.
\[
|\Htp(s,x,y)| \leq |\Htp(s,x,y)|^{\frac 12} |\Htp(s,x,y)|^{\frac 12}
\leq \frac{C}s  D(s,x,y).
\]

Assume the result holds for $H_{n-1}$. 
Let 
$r = \frac{\sqrt{s_0}} {16}$.
Let $\psi\in C^\infty_c\big(Q_{2r}(s_0,y_0)\big)$
where $\psi\Big|_{Q_r(s_0,y_0)}\equiv 1$, $0\leq\psi\leq 1$, and
$\frac {\p^j\psi}{\p s^j} \leq \frac{c_j}{r^{2j}}$. We can use Lemma \ref{lem:kur sub} because
if $s>0$, $H_{n-1}(s,z,w)$ satisfies $\H H_{n-1}(s,x,y)=0$.
Using Lemma \ref{lem:kur sub} and Proposition \ref{prop:L^2 s deriv}, for
$r>0$ and $Q = Q_{2r}(s_0,y_0)$
\begin{align*}
&\left|\frac{\p^n \Htp}{\p s^n}(s_0,x,y_0)\right| 
\leq \frac C{r^2}\left( \iint_{Q_r(s_0,y_0)}\left|\frac{\p^n \Htp}{\p s^n}(s,x,y)\right|^2\, ds dy\right)^{\frac 12}\\
&\leq \frac C{r^2}\left( \iint_{\R\times\C} \frac{\p^n \Htp}{\p s^n}(s,x,y)
\overline{\frac{\p^n\Htp}{\p s^n}(s,x,y)} \psi(s,y)\, ds dy\right)^{\frac 12}\\
&=\frac{C}{r^2}\left( \iint_{\R\times\C} \overline{\Htp(s,x,y)} \sum_{j=0}^{n}
\frac{\p^{n+j} \Htp}{\p s^{n+j}}(s,x,y) \frac{\p^{n-j}\psi}{\p s^{n-j}}(s,y)\, ds dy\right)^{\frac 12} \\
&\leq \frac{C}{r^2}\left[ \|\Htp(\cdot,x,\cdot)\|_{L^2(Q)}\sum_{j=0}^{n} c_j \frac{1}{r^{2(n-j)}}
\left\|\frac{\p^{n+j} \Htp}{\p s^{n+j}}(\cdot,x,\cdot)\right\|_{L^2(Q)}\right]^{\frac 12} \\
&\leq \frac{C}{r^2}\left[\Htp(s_0,x,y_0) r^2\left(\frac{1}{s_0^{2n} }
+ \frac{1}{r^{2n} s_0^n}\right)\right]^{1/2}\\
&\leq \frac{C_n}{r} \frac{D(s_0,x,y_0)^{\frac12}}{s_0^{\frac 12}}
\left(\frac{1}{s_0^{n} } + \frac{1}{r^{n} s_0^{\frac n2}}\right)
\leq \frac{C_n}{s_0^{n + 1}} D(s_0,x,y_0)^\frac{1}{2}.
\end{align*}
\end{proof}

\n Integrating in $y$ gives the immediate corollary:
%Corollary: L^2 estimates on s derivatives of H
\begin{cor}\label{cor:L^2 estimates on s derivatives}
Let $n_1$, $n_2$, $n_3\geq 0$ and $n=n_1+n_2+n_3$. Then there exists $C_n>0$ so that
\[
\left\|\frac{\p^{n_1}}{\p s^{n_1}} \Boxtpx^{n_2} (\Boxtpy\sh)^{n_3} \Htp(s,x,\cdot)\right\|_{L^2(\C)} 
\leq \frac{C_n}{s^{n+\frac12}}. 
\]
\end{cor}
\qed
%Lemma: X^\alpha H(s,x,y)
\begin{lem}\label{lem:X^alpha H} Let $\alpha$ be a multiindex and $j\geq 0$. Then there exists $C_{|\alpha|,j}>0$
so that if $R = \min\{\frac{\sqrt{s_0}}{16},\frac{\mu(x_0,\tfrac 1\tau)}4\}$, then {\small
\[
\big| X_x^\alpha(\Boxtpy\sh)^j\Htp(s,x,y)\big| + \big|\Boxtpx^j  U_y^\alpha\Htp(s,x,y)\big|
\leq \frac{C_{|\alpha|}}{R^{\frac 12+\frac12|\alpha|} s^{\frac 34+j+\frac 14|\alpha|}} D(s,x_0,y)^{\frac 14}.
\] }
\end{lem}

\begin{proof} It is enough to bound $|U_y^\alpha \Boxtpx^j \Htp(s,x_0,y)|$
for a fixed $x_0\in\C$.
In fact,
we can even assume that $\frac{\p^n p}{\p z^n}(x_0) = \frac{\p^n p}{\p\z^n}(x_0)=0$ for all $n$ by
Proposition \ref{prop:rotation}.
This means if $|y-x_0|\leq \mu(x_0,1/\tau)$,
\begin{equation}\label{eq:pjk bound}
\left|\frac{\p^{j+k}p}{\p z^j \p\z^k}(y)\right| \leq \frac{1}{\mu(x_0,1/\tau)^{j+k}}\tau \Lambda(x_0,\mu(x_0,1/\tau))
\sim \frac{1}{\mu(x_0,1/\tau)^{j+k}}.
\end{equation}
Let $R = \min\{\frac{\sqrt s}{16}, \frac 14\mu(x_0,\frac1\tau)\}$. Also, fix $s$ and let
$g(y) = \Boxtpx^j \Htp(s,x_0,y)$. Let $D$ stand for $\frac{\p}{\p x}$ or $\frac{\p}{\p y}$.
Then from Theorem \ref{thm:sobolev 1}, if $\vp \in C^\infty_c(D(x_0,R))$ 
with $0 \leq \vp \leq 1$,
$|D^\beta\vp| \leq \frac {c_{|\beta|}}{R^{|\beta|}}$ for $0 \leq |\beta| \leq 2$, we have
\begin{equation}\label{eqn:sobog}
|U_y^\alpha g(y)| \leq \frac{C}{R} \sum_{|\beta|\leq 2} R^{|\beta|} \|\vp^{\frac 12}
U_y^\beta U_y^\alpha g\|_{L^2(\C)}.
\end{equation}
Then
\begin{align}
 \|\vp^{\frac 12}&U_y^\beta U_y^\alpha g\|_{L^2(\C)}^2
\leq \Big( U_y^\beta U_y^\alpha g, \vp U_y^\beta U_y^\alpha g\Big)
=\left| \Big(g, U_y^\beta U_y^\alpha\big(\vp U_y^\beta U_y^\alpha g\big)\Big)\right|\nn\\
%&= \sum_{|\gamma_1|+|\gamma_2| = |\alpha|+|\beta|} c_{\gamma_1,\gamma_2} 
%\Big(g, (D^{\gamma_1}\vp)U_y^{\gamma_2}U_y^\beta U_y^\alpha g\Big)\nonumber\\
&\leq \sum_{|\gamma_1|+|\gamma_2| = |\alpha|+|\beta|} c_{\gamma_1,\gamma_2} 
\|g\|_{L^2(D(x_0,R))} \frac{1}{R^{|\gamma_1|}}
\|U_y^{\gamma_2}U_y^\beta U_y^\alpha g\|_{L^2(D(x_0,R))}.\label{eqn:L^2g}
\end{align}
Next, from Corollary \ref{cor:L^2 estimates on s derivatives}, Proposition \ref{prop:X to B Box}, and
Theorem \ref{thm:semigroups}, for some order 0 OPF operator $B_\tau$,
we have the estimate (note the complex conjugate in the first inequality),
\begin{align}
\|U_y^{\gamma_2}U_y^\beta &U_y^\alpha g\|_{L^2(D(x_0,R))}^2
\leq \Big(  X_y^{\gamma_2} X_y^\beta X_y^\alpha \bar g, 
 X_y^{\gamma_2} X_y^\beta X_y^\alpha \bar g \Big)\nn\\
&= \Big|\Big( \bar g,  X_y^\alpha X_y^\beta X_y^{\gamma_2}
 X_y^{\gamma_2} X_y^\beta X_y^\alpha \bar g \Big)\Big|\nonumber \\
&\leq \|g\|_{L^2(\C)}\|B_\tau \Boxtp^{|\gamma_2|+|\beta|+|\alpha|} \bar g\|_{L^2(\C)}
\leq \frac{C}{s^{\frac 12+ j}} s^{-(|\gamma_2|+|\beta|+|\alpha| + \frac 12)}.\label{eqn:L^2g2}
\end{align}
Plugging \eqref{eqn:L^2g} into \eqref{eqn:L^2g2}  gives
\begin{equation} \label{eqn:L^2g3}
\|\vp^{\frac 12} X_y^\beta X_y^\alpha \bar g\|_{L^2(\C)}^2
\leq C |g(y)|R\sum_{|\gamma_1|+|\gamma_2| = |\alpha|+|\beta|} c_{\gamma_1,\gamma_2} 
\frac{1}{R^{|\gamma_1|}}s^{-\frac 12(|\gamma_2|+|\beta|+|\alpha| + 1+2j)}.
\end{equation}
Using the fact that $R \leq \sqrt s$ and inserting \eqref{eqn:L^2g3} into \eqref{eqn:sobog},
we have
\begin{align*}
| X_y^\alpha \overline{g(y)}| 
&\leq  \frac {C_{|\alpha|}}{R}\sum_{|\beta|\leq 2} R^{|\beta|} |g(y)|^{\frac12}R^{\frac 12}
\sum_{|\gamma_1|+|\gamma_2| = |\alpha|+|\beta|} s^{-\frac 14(|\gamma_2|+|\beta|+|\alpha| + 1+2j)}
\frac{1}{R^{\frac{|\gamma_1|}2}}\\
&\leq \frac {C_{|\alpha|}}{R^{\frac 12} s^{\frac 34+j}} D(s,x_0,y)^{\frac 14} \sum_{|\beta|\leq 2}
\sum_{|\gamma_1|+|\gamma_2| = |\alpha|+|\beta|} R^{|\beta|-\frac 12|\gamma_1|}s^{-\frac 14|\gamma_2|}
s^{-\frac14|\beta|}s^{-\frac 14|\alpha|}\\
&\leq \frac {C_{|\alpha|}}{R^{\frac 12} s^{\frac 34+j}} D(s,x_0,y)^{\frac 14} R^{-\frac12|\alpha|}
s^{-\frac 14|\alpha|}.
\end{align*}
\end{proof} 

%Corollary: L^2 estimates based on X^\alpha\Box^j estimates
\begin{cor}\label{cor:L^2 X^alpha Box^j estimates}
Let $\alpha$ be a multiindex and $j\geq 0$. Then there exists $C_{|\alpha|,j}>0$
so that
\[
\big\|  X_x^\alpha (\Boxtpy\sh)^j\Htp(s,x,\cdot)\big\|_{L^2(\C)} 
+ \big\|U_y^\alpha\Boxtpx^j \Htp(s,x,\cdot)\big\|_{L^2(\C)}
\leq \frac{C_{|\alpha|,j}}{s^{\frac 12+j+\frac{|\alpha|}2}}.
\]
\end{cor}

\begin{proof} Using the estimate from Lemma \ref{lem:X^alpha H},
if $R = \frac{\sqrt{s}}{16}$, then the result follows by direct calculation and a simple change of
variables. If $R = \frac 14\mu(x,1/\tau)$, then we use the fact that 
$D(s,x,y)^{\frac 14} \leq C_j D(s,x,y)^{\frac 18}\left( \frac{\mu(x,1/\tau)^2}{s}\right)^j$  for any $j\geq 0$.
With this estimate, the result follows immediately.
\end{proof}

The final lemma we need is:
%Lemma: all derivatives
\begin{lem}\label{lem:X^alpha X^beta H}
Let $\alpha$ and $\beta$ be multiindices. 
If $R = \min\{\frac{\sqrt{s_0}}{16},\frac{\mu(x_0,1/\tau)}4\}$,
there exists $C_{|\alpha|,|\beta|}>0$ so that
\[
\big|  X_x^\alpha  X_y^\beta \Htp(s,x,y)\big| 
\leq C_{|\alpha|,|\beta|}\frac{1}{R^{\frac 34}s^{\frac 58}} R^{-\frac{|\alpha|}2 - \frac{|\beta|}4}
s^{-\frac{|\alpha|}4 - \frac{3|\beta|}8}D(s,x,y)^{\frac 18}.
\]
\end{lem}

\begin{proof} As in Lemma \ref{lem:X^alpha H}, we may assume that 
$\frac{\p^n p}{\p z^n}(x_0) = \frac{\p^n p}{\p\z^n}(x_0)=0$ for all $n$ so by
\eqref{eq:pjk bound} 
\[
\left|\frac{\p^{j+k}p}{\p z^j \p\z^k}(y)\right| \les \frac{1}{\mu(x_0,1/\tau)^{j+k}}.
\]
Fix $s$ and $x_0$. Let $y\in D(x_0,R)$. Let $\vp\in C^\infty_c(D(x_0,2R))$ so that 
$\vp\Big|_{D(x_0,R)}\equiv 1$, $0\leq \vp \leq 1$ and $|D^{\alpha}\vp|\leq 
\frac{c_{|\alpha|}}{R^{|\alpha|}}$. Let $f(x) =  X_y^\beta \Htp(s,x,y)$ and 
$g(x) =  X_x^\alpha  X_y^\beta \Htp(s,x,y)$. From Theorem \ref{thm:sobolev 1},
\[
|g(x_0)| \leq \frac{C}{R}\sum_{|\gamma|\leq 2} R^{|\gamma|}
\|\vp^{\frac 12}  X_x^\gamma g\|_{L^2(\C)}.
\]
Next,
\begin{align*}
\|\vp^{\frac 12}  X_x^\gamma g\|_{L^2(\C)}^2
&= \Big(  X_x^\gamma X_x^\alpha f, \vp  X_x^\gamma g \Big)
= \left|\Big(f, X_x^\alpha X_x^\gamma\big[\vp  X_x^\gamma g\big]\Big)\right|\\
&= \sum_{|\gamma_1|+|\gamma_2| = |\gamma|+|\alpha|}c_{\gamma_1,\gamma_2}
\left|\Big(f, D^{\gamma_1}\vp  X_x^{\gamma_2} X_x^\gamma g\Big) \right| \\
&\leq \sum_{|\gamma_1|+|\gamma_2| = |\gamma|+|\alpha|}c_{\gamma_1,\gamma_2}
\| f \|_{L^2(D(x_0,R))} \frac{1}{R^{|\gamma_1|}}
\| X_x^{\gamma_2} X_x^\gamma g\|_{L^2(\C)}.
\end{align*}
Using Proposition \ref{prop:X to B Box} and Corollary \ref{cor:L^2 X^alpha Box^j estimates}, for
some order zero OPF operator $B_\tau$ we have
\begin{align*}
\| X_x^{\gamma_2}& X_x^\gamma g\|_{L^2(\C)}^2
= \Big(  X_x^{\gamma_2} X_x^\gamma X_x^{\alpha}f, 
 X_x^{\gamma_2} X_x^\gamma X_x^{\alpha}f\Big) 
= \left| \Big(f,  X_x^{\alpha} X_x^\gamma X_x^{\gamma_2}
 X_x^{\gamma_2} X_x^\gamma X_x^{\alpha}f\Big)\right| \\
&\leq \|f\|_{L^2(\C)} \|B_\tau \Boxtp^{|\alpha|+|\gamma|+|\gamma_2|}f\|_{L^2(\C)}
%\leq C_{|\alpha|+|\gamma|+|\gamma_2|+|\beta|} s^{-\frac 12(|\beta|+1)}
%s^{-\frac12(|\beta|+1)-|\alpha|-|\gamma|-|\gamma_2|}\\
= C_{|\alpha|+|\gamma|+|\gamma_2|+|\beta|} 
s^{-|\beta|-1-|\alpha|-|\gamma|-|\gamma_2|}.
\end{align*}
Thus, since $\|f\|_{L^2(D(x_0,R))} \leq C |f(x_0)|R$,
 {\small\begin{align*}
&|g(x_0)|\leq C_{|\alpha|,|\beta|}\frac1R \sum_{|\gamma|\leq 2} 
\sum_{|\gamma_1|+|\gamma_2| = |\gamma|+|\alpha|} R^{|\gamma|} R^{-\frac{|\gamma_1|}2}
s^{-\frac 14(|\beta|+1+|\alpha|+|\gamma|+|\gamma_2|)}|f(x_0)|^{\frac 12} R^{\frac 12} \\
&\leq C_{|\alpha|,|\beta|}\frac1{R^{\frac 34}s^{\frac 58}} \sum_{|\gamma|\leq 2} 
\sum_{|\gamma_1|+|\gamma_2| = |\gamma|+|\alpha|} 
\hspace{-23.74045pt}R^{|\gamma|} R^{-\frac{|\gamma_1|}2}
s^{-\frac{|\gamma_2|}4} s^{-\frac 14(|\gamma|+|\alpha|)} s^{-\frac{|\beta|}4}
s^{-\frac{|\beta|}8}R^{-\frac{|\beta|}4}D(s,x,y)^{\frac 18}\\
&\leq C_{|\alpha|,|\beta|}\frac1{R^{\frac 34}s^{\frac 58}} R^{-\frac{|\alpha|}2}s^{-\frac{|\alpha|}4}
R^{-\frac{|\beta|}4}s^{-\frac{3}8|\beta|}D(s,x,y)^{\frac 18}.
\end{align*} }
\end{proof}

As an immediate consequence of Lemma \ref{lem:X^alpha X^beta H}, we have:
%
% THEOREM: pointwise bounds on the derivatives of H
%
\va
\n\textbf{Theorem \ref{thm:bounds on H}}
\emph{Let $p$ be a subharmonic, nonharmonic polynomial and $\tau>0$ a parameter.
If $n\geq0$ and $Y^\alpha$ is a product of $|\alpha|$ operators 
$Y = \Zbstp$ or $\Zstp$ when acting in $z$ and $\overline{(\Zstp)}$ or $\overline{(\Zbstp)}$ 
when acting in $w$,
there exist constants $c,c_1>0$ independent of $\tau$ so that
\[
\left|\frac{\p^n}{\p s^n}Y^\alpha  H_{\tau p}(s,z,w)\right|
\leq c_1 \frac{1}{s^{n + \frac 12|\alpha|+1}}
e^{-\frac{|z-w|^2}{32s}} e^{-c \frac{s}{\mu(z,1/\tau)^2}}
e^{-c \frac s{\mu(w,1/\tau)^2}}.
\]
Also, $c$ can be taken with no dependence on $n$ and $\alpha$.}

\va

\begin{proof} The theorem follows Lemma \ref{lem:X^alpha X^beta H}
using the argument of the proof of Corollary \ref{cor:L^2 X^alpha Box^j estimates}.
Also, by the argument of Lemma \ref{lem:X^alpha X^beta H}, specifically 
the argument of (\ref{eqn:sobog})-(\ref{eqn:L^2g2}), we may take $c$
to be independent of $\alpha$ and $n$.
\end{proof}

Using Theorem \ref{thm:bounds on H}, we can integrate in $s$ and recover estimates on 
$G_{\tau p}(z,w)$, the fundamental solution of $\Boxtp$. 
% Corollary: estimates on G(z,w)
\va
\n\textbf{Corollary \ref{cor:bounds on Box inv}} 
\emph {Let $G_{\tau p}(z,w)$ be the integral kernel of the fundamental solution
for $\Boxtp^{-1}$. If $X^\alpha$ is a product of $|\alpha|$ operators of the form
$X^j = \Zbstp, \Zstp$ if acting in $z$ and $\overline{(\Zbstp)}, \overline{(\Zstp)}$ if acting in $w$, then
there exists constants $C_{1,|\alpha|}, C_{2}>0$
so that if $\tau>0$,
\[
|X^\alpha G_{\tau p}(z,w)| \leq C_{1,|\alpha|} \begin{cases}
\log\Big(\frac{2\mu(z,1/\tau)}{|z-w|}\Big)   & |z-w| \leq \mu(z,1/\tau),\ |\alpha|=0\vspace*{6pt} \\
|z-w|^{-|\alpha|} & |z-w| \leq \mu(z,1/\tau),\ |\alpha|\geq 1\vspace*{6pt} \\
\displaystyle\frac{e^{-C_2 \frac{|z-w|}{\mu(z,1/\tau)}}e^{-C_2 \frac{|z-w|}{\mu(w,1/\tau)}}}
{\mu(z,1/\tau)^{|\alpha|}} 
 & |z-w| \geq \mu(z, 1/\tau). \end{cases}
\]
Also, $C_2$ does not depend on $\alpha$.}

\va

\begin{proof} We just need to integrate in $s$ for the estimate. We first show the $|\alpha|=0$ case.
Let $\delta>0$. Then if $c_2 = \frac 1{32}$, 
\[
\int_{0}^\infty \Htp(s,x,y)\, ds
\leq \int_{0}^{\delta} \frac 1s e^{-c_2\frac{|x-y|^2}{s}}\, ds
+ \int_\delta^\infty \frac 1s e^{-c \frac{s}{\mu(x,1/\tau)^2}} \, ds= I + II.
\]
To estimate $I$, we let $t = c_2 \frac{|x-y|^2}s$, so $-\frac 1t\, dt = \frac 1s\, ds$ and
\[
I = \int_{\frac{|x-y|^2}\delta}^\infty \frac 1t e^{-t}\, dt.
\]
If $c_2\frac{|x-y|^2}\delta \leq 1$, then
\[
I = \int_{c_2\frac{|x-y|^2}\delta}^1 \frac 1t e^{-t}\, dt + \int_{1}^\infty \frac 1t e^{-t}\, dt
\leq C\left (\log\Big(\frac\delta{|x-y|^2}\Big) + 1\right).
\]
Also, if $c_2\frac{|x-y|^2}\delta\geq 1$,
\[
I \leq \frac{1}{c_2\frac{|x-y|^2}\delta} \int_{c_2\frac{|x-y|^2}\delta}^\infty e^{-t}\, dt
= C\frac\delta{|x-y|^2} e^{-c_2\frac{|x-y|^2}\delta} \leq Ce^{-c_2\frac{|x-y|^2}\delta}.
\]
To estimate $II$, set $t= c \frac{s}{\mu(x,1/\tau)^2}$, and we have
\[
II = \int_{c\frac{\delta}{\mu(x,1/\tau)^2}}^\infty \frac 1t e^{- t}\, dt.
\]
If $c\frac{\delta}{\mu(x,1/\tau)^2} \leq 1$, we have
\[
II = \int_{c\frac{\delta}{\mu(x,1/\tau)^2}}^1 \frac 1t e^{-t}\, dt 
+ \int_{1}^\infty \frac 1t e^{-t}\, dt
\leq C\left (\log\Big(\frac{\mu(x,1/\tau)^2}\delta\Big) + 1\right).
\]
Also, if $c\frac\delta{\mu(x,1/\tau)^2}\geq 1$,
\[
II \leq \Big(c\frac\delta{\mu(x,1/\tau)^2}\Big)^{-1}
 \int_{c\frac{\delta}{\mu(x,1/\tau)^2}}^\infty \hspace{-11.14pt}e^{-t}\, dt
= \frac{\mu(x,1/\tau)^2}\delta e^{-c\frac\delta{\mu(x,1/\tau)^2}}
\leq C e^{-c\frac\delta{\mu(x,1/\tau)^2}}.
\]
Setting $\delta = \frac{|x-y|}{\mu(x,1/\tau)}$ yields the result. The $|\alpha|\geq 1$ case uses
the same argument as the $|\alpha|=0$ case, 
except that the on-diagonal estimate is simpler since $\int_{\delta}^\infty 
\frac{1}{s^{1+\frac 12|\alpha|}}\, ds \sim s^{-\frac 12|\alpha|}$ converges.
\end{proof}

%%%%%%%%%%%%5
%
%
%  SECTION: COMPARISON OF MY ESTIMATES WITH CHRIST'S
%
%
%%%%%%%%%%%%%5
\appendix
\section{A Comparison of the Estimates of $G_{p}(z,w)$}\label{sec:est equiv}

We will focus on homogenous polynomials of the form $p_1(z) = |z|^{2m}$ and $p_2(x+iy) = x^{2m}$. 
We will write $z=x+iy$. Recall from above that if $\rho_{p_j}$ is the metric
in \cite{Christ91} associated to $p_j$, the $d\rho_{p_j}^2 \sim \mu_{p_j}(z,1)^{-2}ds^2$.

To show that the estimates in Corollary \ref{cor:bounds on Box inv} 
and Christ's estimates in \cite{Christ91}
of $\G_{p}(z,w)$ agree, we must show that for $j=1$ or $2$,
\begin{equation}\label{eq:est equiv}
\frac{|z-w|}{\mu_{p_j}(w,1)} + \frac{|z-w|}{\mu_{p_j}(z,1)} \sim \rho_{p_j}(z,w).
\end{equation}
\eqref{eq:est equiv} will follow from
Corollary \ref{cor:my G est size} and Proposition \ref{prop:rho} for both
$p_1$ and $p_2$.
In fact, the $p_1$ case shows that the estimates agree 
whenever $p(z)$ is a homogenous, subharmonic polynomial
of degree $2m$ whose Laplacian does not vanish on the unit circle. Also, in the case $m=2$, an 
elementary computation shows that
$p(z) = x^4$ is equivalent to 
the general case for subharmonic, nonharmonic homogenous polynomials of degree 4.

\begin{prop}\label{prop:l_j} 
\[
 \mu_{p_1}(z,1)\sim \min\{1,\frac{1}{|z|^{m-1}}\},
\qquad \mu_{p_2}(z,1)\sim \min\{1,\frac{1}{|x|^{m-1}}\}.
\]
\end{prop}

This computation has the following immediate corollaries.
\begin{cor}\label{cor:my G est size} 
\begin{align*}
\frac{|z-w|}{\mu_{p_1}(w,1)} + \frac{|z-w|}{\mu_{p_1}(z,1)} &\sim |z-w| + |z-w|(|z|^{m-1} + |w|^{m-1})\\
\intertext{and}
\frac{|z-w|}{\mu_{p_2}(w,1)} + \frac{|z-w|}{\mu_{p_2}(z,1)} &\sim |z-w| + |z-w|(|\Rre z|^{m-1} + |\Rre w|^{m-1}).
\end{align*}
\end{cor}

\begin{cor}\label{prop:lj sum} 
\[
 \frac{1}{\mu_{p_1}(z,1)^2} \sim 1+|z|^{2m-2}, \qquad \frac{1}{\mu_{p_2}(z,1)^2} \sim 1+|x|^{2m-2}.
\]
\end{cor}

At this point, we will concentrate on sketching a computation of 
$\rho_{p_1}(z,w)$. The computation for  $\rho_{p_2(z,w)}$ is analogous and shows
$\rho_{p_2}(z,\zeta) \sim |z-\zeta|+ |z-\zeta|(|\Rre z|^{m-1}+|\Rre \zeta|^{m-1})$.

\begin{prop}\label{prop:rho}
\[
 \rho_{p_1}(z,\zeta) \sim |z-\zeta|+ |z-\zeta|(|z|^{m-1}+|\zeta|^{m-1}),
\]
where the constant depends only on $m$.
\end{prop}

\begin{proof} As a consequence of Corollary \ref{prop:lj sum},
\begin{equation}\label{equ:rho equiv}
\rho_{p_1}(z,\zeta) \sim \inf_{\alpha}\left\{\int_{0}^{1}
(1+|\alpha(t)|^{m-1})|\alpha'(t)|\,dt\right\}
\end{equation}
If $\alpha(t) = z(1-t) + \zeta t$, then
\[
\int_{0}^{1} (1+|\alpha(t)|^{m-1})|\alpha'(t)|\,dt
\les |z-\zeta|+ |z-\zeta|(|z|^{m-1}+|\zeta|^{m-1}).
\]
For the other direction, we give a more complete argument.
{\small
\[
 \rho(z,\zeta) \sim \inf_{\alpha}\left\{\int_{0}^{1}
  (1+|\alpha(t)|^{m-1})|\alpha'(t)|\,dt\right\}
 \ges |z-\zeta| + \inf_{\alpha}\left\{\int_{0}^{1}
  |\alpha(t)|^{m-1}|\alpha'(t)|\,dt\right\}
\] }
Now, set $\gamma(z,\zeta) = \inf_{\alpha}\left\{\int_{0}^{1}
|\alpha(t)|^{m-1}|\alpha'(t)|\,dt\right\}$. Then
\[
 \gamma(z,\zeta)
 = \inf_{\alpha}\left\{\int_{0}^{1}\left|\frac{d}{dt}
  \alpha(t)^{m}\right|\,dt\right\}
  \geq \inf_{\atopp{\beta}{\atopp{\beta(0)=z^{m}}{\beta(1)=\zeta^{m}}}}
 \left\{\int_{0}^{1}|\beta'(t)|\,dt\right\} = |z^{m}-w^{m}|.
\]
Thus, if $|z|\geq 2|\zeta|$, then $|z^{m}-\zeta^{m}|\sim |z^{m}|
\sim|z-\zeta|(|z|^{m-1}+|\zeta|^{m-1})$. Also, note that there is a 
one-to-one correspondence of paths between $z$ and $\zeta$ and paths
between $rz$ and $r\zeta$ by sending $\alpha(t)$ to $\alpha_{r}(t) = 
r\alpha(t)$. Hence, it follows immediately that $\gamma(rz,r\zeta)
= r^{m}\gamma(z,\zeta)$, so without loss of generality, we can assume
that $|z|=1$ and $\frac 12 <|\zeta| \leq 1$. 
Let $\frac \zeta z = re^{it}$. Note that
$t = \arg \zeta-\arg z$. Also, we can write
\[
 (z^m-\zeta^m) = z^m(1-r^m e^{imt}) 
 = z^m(1-r e^{it})\prod_{k=1}^{m-1}(e^{2\pi i\frac{k}{m}} - re^{it}).
\]
Note that if $|t|<\frac{\pi}m$, then 
$|e^{2\pi i\frac{k}{m}} - re^{it}| >c>0$ for some constant $c$
and for all $k=1,2,\ldots,m-1$, and in this
case
\[
 |z^m(1-r e^{it})\prod_{k=1}^{m-1}(e^{2\pi i\frac{k}{m}} - re^{it})|
 \geq c |z-\zeta| |z|^{m-1}\geq \frac c2 |z-\zeta|(|z|^{m-1} + 
 |\zeta|^{m-1}).
\]
Finally, if $|t|\geq \frac{\pi}m$, then $\rho(z,\zeta)\geq c$ for some
constant $c>0$, and this is the desired result as $|z-\zeta|\geq c_1$ and
$|z|,|\zeta| \in [\frac 12,1]$.
\end{proof}
% \appendix

% \section{}
% \label{}

% Bibliographic references with the natbib package:
% Parenthetical: \citep{Bai92} produces (Bailyn 1992).
% Textual: \citet{Bai95} produces Bailyn et al. (1995).
% An affix and part of a reference:
%   \citep[e.g.][Ch. 2]{Bar76}
%   produces (e.g. Barnes et al. 1976, Ch. 2).
\bibliographystyle{plain}
%\bibliography{mybib}

%\begin{thebibliography}{}

%% \bibitem[Names(Year)]{label} or \bibitem[Names(Year)Long names]{label}.
%% (\harvarditem{Name}{Year}{label} is also supported.)

%\bibitem[]{}

%\end{thebibliography}

\end{document}